\documentclass[a4paper, 11pt]{article}
\usepackage{amssymb,amsmath,graphicx}
\usepackage{amsfonts}
\usepackage{epsfig,latexsym,graphicx}
\newtheorem{theorem}{Theorem}[section]

\newtheorem{lemma}[theorem]{Lemma}
\newtheorem{proposition}[theorem]{Proposition}

\newtheorem{definition}{Definition}[section]

\begin{document}

\title{\bf{Robust and Efficient Parameter Estimation based on Censored Data with Stochastic Covariates
}
}
\author{Abhik Ghosh and Ayanendranath Basu \\
Indian Statistical Institute \\
{\it abhianik@gmail.com, ayanbasu@isical.ac.in}}
\date{}
\maketitle

\begin{abstract}
Analysis of random censored life-time data along with some related stochastic covariables
is of great importance in many  applied sciences like medical research, population studies and planning etc.
The parametric estimation technique commonly used under this set-up is based on the efficient but non-robust 
likelihood approach. In this paper, we propose a robust parametric estimator 
for the censored data with stochastic covariates based on the minimum density power divergence approach.
The resulting estimator also has competitive efficiency with respect to the maximum likelihood estimator 
under pure data. The strong robustness property of the proposed estimator with respect to the presence of 
outliers is examined and illustrated  through an appropriate simulation study in the context of 
censored regression with stochastic covariates. 
Further, the theoretical asymptotic properties of the proposed estimator are also derived 
in terms of a general class of M-estimators based on the estimating equation.
\end{abstract}

\textbf{Keywords:}  Censored Data; Robust Methods; Linear Regression; Density power divergence; M-Estimator;
Exponential Regression Model, Accelerated Failure Time Model.

\section{Introduction}\label{SEC:intro}

It is often necessary  to analyze life-time data in many applied sciences including medical sciences, 
population studies, planning etc. For these survival analyses problems, researchers often cannot observe
the full data because some of the respondents may leave the study in between or 
some may be still alive at the end of the study period. Statistical modelling of such data involves the 
idea of censored distributions and random censoring variables. Mathematically,
let $Y_1, \ldots, Y_n$ be $n$ independent and identically distributed (i.i.d.) observations from 
the population with unknown life-time distribution $G_{Y}$.
We assume that the observations are censored by a censoring distribution $G_C$ independent of $G_Y$ 
and $C_1, \ldots, C_n$ denote $n$ i.i.d.~sample observations from $G_C$.
We only observe the portion of $Y_i$s (right) censored by $C_i$s, i.e., we observe
$$
Z_i = \min \left(Y_i, C_i\right) ~~~ \mbox{and}~~~ \delta_i = I(Y_i \leq C_i), ~~ i=1, \ldots, n,
$$
where $I(A)$ denote the indicator function of the event $A$. Based on these data $(Z_i, \delta_i)$,
our aim is to do inference about the lifetime distribution $G_Y$.
Suppose $Z_{(i, n)}$ denotes the $i$-th order statistic in $\{Z_1, \cdots, Z_n\}$ and 
$\delta_{[i, n]}$ is the value of corresponding $\delta$ ($i$-th concomitant).
The famous product-limit (non-parametric) estimator of $G_Y$ under this set-up had been derived by 
Kaplan and Meier (1958), which is given by
$$
\widehat{G_Y}(y) = 1 - \prod_{i=1}^{n}\left[1 - \frac{\delta_{[i, n]}}{n-i+1}\right]^{I(Z_{(i,n)}\leq y)}.
$$
It can be seen that, under suitable assumptions,  the above product-limit estimator is in fact 
the maximum likelihood estimator of the distribution function in presence of censoring
and enjoys several optimum properties. Many researchers have proved such properties and 
also extended it for different complicated inference problems with censored data; 
for example, see Petersen (1977), Chen et al.~(1982), Campbell and F\"{o}ldes (1982), 
Wang et al.~(1986), Tsi et al.~(1987),
Dabrowska (1988), Lo et al.~(1989), Zhou (1991), Cai (1998), Satten and Datta (2001) among many others.

In this present paper we further assume the availability of a set of uncensored covariables 
$ X \in \mathbb{R}^p$ that are associated with our target response $Y$; i.e., 
for each respondent $i$ we have observed the values $X_i$ along with $(Z_i, \delta_i)$.
These covariables are generally the demographic conditions of the subject or some measurable indicator of 
the response variable (e.g., medical diagnostic measures like blood pressure, hemoglobin content, etc., for 
clinical trial responses). Let us assume that the distribution function of 
these i.i.d.~covariates is $G_X$ and their joint distribution function with $Y$ is $G$ so that 
$$
G_Y(y) = \int G(x,y)dx = \int G_{Y|X=x}(y)G_X(x)dx,
$$
where $G_{Y|X=x}$ is the conditional distribution function of $Y$ given $X=x$.
Instead of inferring about the response $Y$ alone, here we are more interested in obtaining the association 
between response and covariables through the conditional distribution $G_{Y|X}$; the distribution $G_X$ of 
covariates is often of interest  but sometimes it may act as a nuisance component too.
Let us denote the $i$-th concomitant of $X$ associated with $Z_{(i,n)}$ by $X_{[i,n]}$.
Under this set-up, Stute (1993) has extended the Kaplan-Meier product limit (KMPL) estimator $\widehat{G_Y}(y)$
to obtain a non-parametric estimate of the joint multivariate distribution function $G$ given by 
$$
\widehat{G}(x,y) = \sum_{i=1}^{n}W_{in} I(X_{[i,n]}\leq x, Z_{(i,n)}\leq y)
$$
where the weights are calculated as 
$$
W_{in} = \frac{\delta_{[i, n]}}{n-i+1} \prod_{j=1}^{i-1}\left[\frac{n-j}{n-j+1}\right]^{\delta_{[j, n]}}.
$$
Note that when there is no censoring at all, i.e., $\delta_i=1$ for all $i$, then $W_{in}=\frac{1}{n}$ for each $i$ 
so that $\widehat{G_Y}(y)$ and $\widehat{G}(x,y) $ coincide with respective empirical distribution functions.
Further, above estimator is also self-adjusted in presence of any ties in the data.
These give the framework for non-parametric inference based on the censored data with covariates. 
Stute (1993, 1996) proved several asymptotic properties like strong consistency, asymptotic distribution 
of  $\widehat{G}(x,y) $ and related statistical functionals.
Several non-parametric and semi-parametric inference procedures using $\widehat{G}(x,y)$
are widely used in real life applications.

However, for many applications in medical sciences, one may know the parametric form of 
the distribution of the survival time (censored responses) possibly through previous experience 
in similar context (similar drugs or similar diseases may have studied in the past). 
In such cases, the use of a fully parametric model is much more appropriate over the semi-parametric 
or non-parametric models. Many advantages of a fully parametric model for the regression with censored responses
had been illustrated in Chapter 8 of Hosmer et al.~(2008) which include -- 
(a) greater efficiency due to the use of full likelihood, 
(b) more meaningful estimates of clinical effects with simple interpretations,
(c) prediction of the response variable from the fitted model etc.
The most common parametric models used for the analysis of survival data are 
the exponential, Weibull or log-logistic distributions.
Many researchers had used such parametric models to analyze survival data more efficiently;
see for example, Cox and Oakes (1984), Crowder et al.~(1991), Collett (2003), Lawless (2003),
Klein and Moeschberger (2003) among others.

The robustness issue with the survival data, on the other hand, has got prominent attention very recently. 
The size and availability of survival data has clearly been growing in recent times in biomedical and other industrial studies  
which often may contain few erroneous observations or outliers and it is very difficult to sort out 
those observations in presence of complicated censoring schemes. Some recent attempts have
been made to obtain robust parametric estimators based on survival data without any covariates. 
For example, Wang (1999) has derived the properties of M-estimators for univariate life-time distributions
and Basu et al.~(2006) have developed a more efficient robust parametric estimator by minimizing the 
density power divergence measure (Basu et al., 1998). These estimators, along with the automatic control for  
the effect of outlying observations, provide a compromise between the most efficient classical parametric estimators 
like maximum likelihood or method of moments and the inefficient non-parametric or semi-parametric approaches
provided there is no significant loss of efficiency under pure data. The present paper extends this idea to develop 
such  robust estimators for the model parameters under censored response with covariates.   
It does not follow directly from the existing literature as we need to change the laws of large number 
and central limit theorem for censored data suitably in the presence of covariates.

It is to be noted that, in this paper we consider a fully parametric model for survival data with covariates,
which is not the same as the usual semi-parametric or nonparametric  regression models like 
the Cox proportional hazard model (Cox, 1972) or the Buckley-James linear regression 
(Buckley and James, 1979; Ritov, 1990).
The latter methods are generally more robust with respect to model misspecification
but less efficient compared to a fully parametric models. 
Further, many recent attempts have been made to develop inference under such semi-parametric models 
that are robust also with respect to outliers in the data; e.g., Zhou (1992), Bednarski (1993), Kosorok et al.~(2004),
Bednarski and Borowicz (2006), Salibian-Barrera and Yohai (2008), Farcomeni and Viviani (2011) etc.
However, no such work has been done to develop robust inference under fully parametric regression 
models with survival response. The few works that closely relate to the proposal of our present paper are 
by Zhou (2010), Locatelli et al.~(2011) and Wang et al.~(2015),  
who have proposed some robust solutions for a particular case of accelerated failure time regression model 
without any parametric assumptions on the stochastic covariates;
the first two papers propose the M-estimators and S-estimators respectively for the semi-parametric AFT models
with parametric location-scale error and the third one extends the M-estimators further to achieve robustness 
with respect to leverage points and simultaneous robust estimation of the error variance.
In Zhou (2010) there is the possibility of dependence between the error and stochastic explanatory variables, but independence
is also allowed in which case it considers the same model set-up as in Locatelli et al.~(2011).  
However, none of these approaches considered fully parametric models with suitable assumptions on the marginal covariate distributions.
The present paper fill this gap in the literature of survival analysis by proposing 
a simple yet general robust estimation criterion with more efficiency under any general parametric model
for both the censored response and the stochastic covariates.
Further, the proposal of the present paper is fully general with respect to model assumptions and 
can be easily extended to the semi-parametric models considered in the existing literature; 
we will show that the existing versions of the M-estimators of Zhou (2010) and Wang et al.~(2015) 
can be considered as a special case of our proposal under the semi-parametric extension.
In that sense, our proposed method in the present paper will give a complete general framework for all possible modeling options of 
censored data with stochastic covariates along with the possibility of 
more efficient inference through fully parametric covariate distributions.

The rest of the paper is organized as follows. 
We start with a brief review of background concepts and results about the  
non-parametric estimator $\widehat{G}(x,y)$ and the minimum density power divergence estimators
in Section \ref{SEC:Background}. Next we consider a general parametric set-up for the censored lifetime data 
with covariates as described above and propose the modified minimum density power divergence estimator
for the present set-up in Section \ref{SEC:MDPDE};  its application in the context of simple linear 
and exponential regression models with censored response and for the general parametric accelerated
failure time model are also described in this section.
In Section \ref{SEC:M-estimators}, we derive theoretical asymptotic properties
for a general class of estimators containing the proposed minimum density power divergence estimator 
under the present set-up; this general class of estimators is indeed a suitable extension of the M-estimators.  
The global nature of our proposal and its generality are discussed in Section \ref{SEC:SemiP_Ex} along with 
the illustration of this extension in the semi-parametric set-up. 
Section \ref{SEC:IF} contains the robustness properties of the proposed MDPDE and the general M-estimators 
examined through the influence function analysis for both fully parametric and semi-parametric set-ups.
The performance of the proposed minimum density power divergence estimator  
in terms of both efficiency and robustness is illustrated 
through some appropriate simulation studies in Section \ref{SEC:Numerical}. 
Some remarks on the choice of the tuning parameter in the proposed estimator 
are presented in Section \ref{SEC:selection_alpha}, while the paper ends with a short concluding remark 
in Section \ref{SEC:Conclusion}.

\section{Preliminary Concepts and Results}
\label{SEC:Background}
\subsection{Asymptotic Properties of $\widehat{G}(x,y)$}
\label{SEC:property_hatG}

One of the main barriers to derive any asymptotic results based on survival data 
was the unavailability of limit theorems, like law of iterated logarithm,  law of large number,
central limit theorem etc., under censorship. This problem has been solved in the recently decades 
mainly through the works of Stute and  Wang; see Stute and Wang (1993), Stute (1995) for 
such limit theorems for the censored data without covariates and 
Stute (1993, 1996) for similar results in presence of covariables.
In this section, we briefly describe some results from the later works with covariates
that will be needed in this paper.

Assume the set-up of life-time variable $Y$ censored by an independent censoring variable $C$ 
as discussed in Section 1.
Denote $Z=\min(Y,C)$; the distribution function of $Z$ is given by $G_Z  = 1 - (1-G_Y)(1-G_C)$. 
In order to have the limiting results for this set-up, we need to make the following basic assumptions:

\begin{itemize}
	\item[(A1)] The life-time variable $Y$ and the censoring variable $C$  are independent and 
	their respective distribution functions $G_Y$ and $G_C$ have no jump in common.
	
	\item[(A2)] The random variable $\delta=I(Y\leq C)$ and $X$ are conditionally independent given $Y,$ 
	i.e. whenever the actual life-time is known the covariates provide no further information on censoring.
	More precisely, $P(Y\leq C | X, Y) = P(Y\leq C | Y).$
\end{itemize}

Now, consider a real valued measurable function $\phi$ from $\mathbb{R}^{p+1}$ to $\mathbb{R}^k$ and 
define
\begin{eqnarray}
	S_n = \sum_{i=1}^{n}W_{in} \phi(X_{[i,n]}, Z_{(i,n)}) = \int \phi(x,y)\widehat{G}(dx, dy).
	\label{EQ:S_n}
\end{eqnarray}
This functional $S_n$ forms the basis of several estimators under this set-up. 
The results presented below describe its strong consistency and distributional convergence;
see Stute (1993, 1996)  for their proofs and details.

\begin{proposition}
	\label{PROP:Consistency}
	[Strong Consistency (Stute, 1993)]\\
	Suppose that $\phi(X,Y)$ is integrable and Assumptions (A1) and (A2) hold for the above mentioned set-up.
	Then we have, with probability one and in the mean, 
	\begin{eqnarray}
		\lim\limits_{n\rightarrow\infty} S_n = \int_{Y<\tau_{G_Z}} \phi(X,Y) dP 
		+ I(\tau_{G_Z} \in A) \int_{Y = \tau_{G_Z}} \phi(X, \tau_{G_Z}) dP,
		\label{EQ:SLLN_surv}
	\end{eqnarray} 
	where $\tau_{G_Z}$ denote the least upper bound for the support of $G_Z$ given by 
	$$
	\tau_{G_z} = inf \{z : G_Z(z)=1\},
	$$
	and $A$ denotes the set of atoms (jumps) of $G_Z$.
\end{proposition}

The above convergence can be written in a simpler form, by defining 
$$
\widetilde{G}(x,y) = \left\{\begin{array}{l c l }
G(x,y) &~\mbox{if}~&~y<\tau_{G_Z}\\
G(x,\tau_{G_Z}-) + I(\tau_{G_Z}\in A) G(x,\{\tau_{G_Z}\}) &~\mbox{if}~&~y\geq\tau_{G_Z}.
\end{array}\right.
$$
Then, the convergence in (\ref{EQ:SLLN_surv}) yields
$$
\lim\limits_{n\rightarrow\infty} \int \phi(x,y)\widehat{G}(dx, dy) = \int \phi(x,y)\widetilde{G}(dx, dy) 
= \widetilde{S}, ~~\mbox{say}.
$$
Further, note that the independence of $Y$ and $C$ gives $\tau_{G_Z} = \min(\tau_{G_Y}, \tau_{G_C})$,
where $\tau_{G_Y}$ and $\tau_{G_C}$ are the least upper bound of the supports of $G_Y$ and $G_C$ respectively.
So, whenever $\tau_{G_Y} < \tau_{G_C}$ or $\tau_{G_C}=\infty$, the modified distribution function 
$\widetilde{G}$ coincides with $G$ and the estimator $S_n$ becomes a strongly consistent estimator 
of its population counterpart $S=\int \phi(x,y){G}(dx, dy)$. 
Further, it follows that the Glivenko-Cantelli type strong uniform convergence 
of $\hat{G}$ to $G$ holds under assumptions (A1) and (A2); see Corollary 1.5 of Stute (1993).

\begin{proposition}
	\label{PROP:CLT}
	[Central Limit Theorem (Stute, 1996, Theorem 1.2)]\\
	Consider the above mentioned set-up with assumption (A2) and suppose that the measurable function 
	$\phi(X,Y)$ satisfies 
	\begin{itemize}
		\item[(A3)] $\int \left[\phi(X,Z)\gamma_0(Z)\delta\right]^2 dP < \infty,$ \\
		\item[(A4)] $\int \left|\phi(X,Z)\right|C^{1/2}(w)\widetilde{G}(dx, dw) < \infty,$
	\end{itemize}
	where 
	\begin{eqnarray}
		\gamma_0(z) &=& \exp\left\{\int_0^{z-}\frac{G_Z^0(dz')}{1-G_Z(z')}\right\}, 
		~~\mbox{with }~~ G_Z^0(z) = P(Z\leq z, \delta=0), \nonumber\\
		\mbox{and }~~ C(w) &=& \int_0^{w-}\frac{G_C(dz')}{[1-G_C(z')][1-G_Z(z')]}.\nonumber
	\end{eqnarray}
	Then we have, as $n\rightarrow\infty$,
	\begin{eqnarray}
		\sqrt{n}(S_n - \widetilde{S}) \mathop{\rightarrow}^\mathcal{D} N(0, \Sigma_\phi),
		\label{EQ:CLT_surv}
	\end{eqnarray} 
	where 
	\begin{eqnarray}
		\Sigma_\phi = Cov\left[\phi(X,Z)\gamma_0(Z)\delta + \gamma_1(Z)(1-\delta) - \gamma_2(Z)\right],
		\label{EQ:CLT_var}
	\end{eqnarray}
	where $\gamma_1$ and $\gamma_2$ are vectors of the same length as $\phi$ and are defined as
	\begin{eqnarray}
		\gamma_1(z) &=& \frac{1}{1-G_Z(z)} \int I(z<w)\phi(x,w)\gamma_0(w)\widetilde{G}^{11}(dx, dw), \nonumber\\
		\mbox{and }~~ \gamma_2(z) &=& \iint \frac{I(v<z,v<w)\phi(x,w)\gamma_0(w)}{[1-G_Z(v)]^2}
		G_Z^0(dv)\widetilde{G}^{11}(dx, dw), \nonumber\\
		\mbox{with }~~ \widetilde{G}^{11}(x, z) &=& P(X\leq x, Z\leq z, \delta=1). \nonumber
	\end{eqnarray}
\end{proposition}

Note that a consistent estimator of the above asymptotic variance can be obtained 
by using the corresponding sample covariance and by replacing the distribution functions 
in the definitions of $\gamma_0$, $\gamma_1$ and $\gamma_2$ by their respective empirical estimators.

In this context, we should note that Assumption (A2) is strictly stronger than the usual
assumptions in regression analysis for censored life-time data (Begun et al., 1983).
This can be seen by writing Stute's estimate $\widehat{G}$ as a particular case of 
the inverse of the probability of censoring weighted (IPCW) statistic, 
where the censoring weights are estimated by the marginal Kaplan-Meier estimator for the censoring time.
This may lead to some biased inference when the censoring distribution depends on the covariates,
but in such cases we cannot have robust results unless moving to the semi-parametric models like Cox regression.
Further, Robins and Rotnitzky (1992) also assumed this stronger condition (A2) 
to develop a more efficient IPCW statistic under the semi-parametric set-up 
(see also Van der Laan and Robins, 2003). 
Zhou (2010) and Wang et al.~(2015) have also considered the same assumption (A2) for robust estimation under semi-parametric 
accelerated failure time models. 
So, while considering the fully parametric set-up throughout this paper, we continue with the assumption (A2)
for deriving any asymptotic result; clearly it does not restrict the practical use of 
the proposed method in finite samples.

\subsection{The Density Power Divergence and Corresponding Estimators}

The density power divergence based statistical inference has become quite popular in recent days 
due to its strong robustness properties and high asymptotic efficiency without using any non-parametric smoothing.
The density power divergence measure between any two densities $g$ and $f$ (with respect to 
some common dominating measure) is defined in terms of a tuning parameter $\alpha \geq 0$ as 
(Basu et al., 1998),
\begin{equation}\label{EQ:dpd}
	d_\alpha(g,f) = \displaystyle \left\{\begin{array}{ll}
		\displaystyle \int  \left[f^{1+\alpha} - \left(1 + \frac{1}{\alpha}\right)  f^\alpha g + 
		\frac{1}{\alpha} g^{1+\alpha}\right], & {\rm for} ~\alpha > 0,\\
		\displaystyle \int g \log(g/f), & {\rm for} ~\alpha = 0.  
	\end{array}\right.
\end{equation}
When we have $n$ i.i.d.~samples $Y_1, \ldots, Y_n$ from a population with true density function $g$, 
modeled by the parametric family of densities $\mathcal{F} = \{f_\theta : \theta \in \Theta \subset \mathbb{R}^p\}$,
the minimum density power divergence estimator (MDPDE) of the parameter $\theta$ 
is to be obtained by minimizing the density power divergence between the data and the model family; 
or equivalently by minimizing 
\begin{eqnarray}
	\int f_\theta^{1+\alpha}(y)dy - \frac{1+\alpha}{\alpha} \int  f_\theta^{\alpha}(y)dG_n(y) 
	=\int f_\theta^{1+\alpha}(y)dy - \frac{1+\alpha}{\alpha} \frac{1}{n}\sum_{i=1}^n  f_\theta^{\alpha}(Y_i),
	\label{EQ:obj_func_MDPDE}
\end{eqnarray}
with respect to $\theta$; here $G_n$ is the empirical distribution function based on the sample.
See Basu et al.~(1998, 2011) for more details and other properties of the MDPDEs.
It is worthwhile to note that the MDPDE corresponding to $\alpha=0$ coincides 
with the maximum likelihood estimator (MLE); the MDPDEs become more robust 
but less efficient as $\alpha$ increases, although the extent of loss is not significant in most cases 
with small positive $\alpha$.  
Thus the parameter $\alpha$ gives a trade-off between robustness and efficiency.
Hong and Kim (2001) and Warwick and Jones (2005) have presented 
some data-driven choices for the selection of optimal tuning parameter $\alpha$.

The MDPDE has been applied to several statistical problems and has been extended suitably 
for different types of data. For example, 
Kim and Lee (2001) have extended it to the case of  robust estimation of extreme value index,
Lee and Song (2006, 2013) have provided extensions in the context of GARCH model and 
diffusion processes respectively
and Ghosh and Basu (2013, 2014) have generalized it to the case of non-identically distributed data  
with applications to the linear regression and the generalized linear model. 
In the context of survival analysis, Basu et al.~(2006) have extended the concept of MDPDE 
for censored data without any covariates to obtain a robust and efficient estimator. 
Based on $n$ i.i.d.~right censored observations $Y_1, \ldots, Y_n$ as above, 
Basu et al.~(2006) have proposed to use the Kaplan-Meier product limit estimator $\widehat{G_Y}$ 
in place of the empirical distribution function $G_n$ in (\ref{EQ:obj_func_MDPDE}) and 
derived the properties of the corresponding MDPDEs.
In the next section, we will further generalize this idea to obtain robust estimators 
for a joint parametric model based on censored data with  covariates.

\section{The Minimum Density Power Divergence Estimation (MDPDE) under Random Censoring with Covariates}
\label{SEC:MDPDE}

\subsection{General Parametric Models and Estimating Equations}\label{SEC:general_MDPDE}

Let us consider the set-up of Section 1. We are interested in making some inference 
about the distribution of the lifetime variable $Y$ and its relation with the covariates 
(through the distribution $G_{Y|X}$) based on the survival data with covariates $(Z_i, \delta_i, X_i)$.
Sometimes one may also be interested in the distribution $G_X$ of the covariates.
As noted earlier, this paper focuses on the parametric approach of inference; 
so we assume two model family of distributions for $G_{Y|X}$ and $G_X$ given by 
$\mathcal{F_X} = \{F_\theta(y|X) : \theta \in \Theta \subseteq \mathbb{R}^q\}$ and 
$\mathcal{F_0} = \{F_{X,\gamma}(x) : \gamma \in \Gamma \subseteq \mathbb{R}^r\}$ respectively.
Then the target parameters of interest are $\theta$ and $\gamma$ which we will estimate jointly based on $(Z_i, \delta_i, X_i)$.
The case of known $\gamma$ can be easily derived from this general case or from the work of Basu et al.~(2006).

The most common and popular method of estimation is the maximum likelihood estimator (MLE)
that is obtained by maximizing the probability of the observed data  $(Z_i, \delta_i, X_i)$
with respect to the parameters $(\theta,~\gamma)$. 
However, in spite of several optimal properties, the MLE has well-known 
drawback of the lack of robustness. As noted in the previous section, the minimum density power divergence 
estimator can be used as a robust alternative to the MLE with no significant loss 
in efficiency under pure data for several common problems.
Deriving the motivation from these, specially from the work of Basu et al.~(2006), 
here we consider the minimum density power divergence estimator (MDPDE) of $(\theta,~\gamma)$ 
obtained by minimizing the objective function (\ref{EQ:obj_func_MDPDE})
for the joint density of the variables $(Y,~X)$ and a suitable estimator of this joint distribution.
Let us denote the density of $F_\theta(y|X)$ by $f_\theta(y|X)$ and so on.
Then the joint model density of $(X,~Y)$ is $f_\theta(y|x)f_{X,\theta}(x)$.
As an estimator of their true joint distribution $G$ we use the KMPL $\widehat{G}(x,y)$, because of its 
optimality properties as described in Section \ref{SEC:property_hatG}. 
Thus, for any $\alpha>0$, the objective function to be minimized with respect to $(\theta,~\gamma)$ is given by
\begin{eqnarray}
	&& H_{n,\alpha}(\theta, \gamma) = \iint f_\theta(y|x)^{1+\alpha}f_{X,\gamma}(x)^{1+\alpha}dxdy 
	- \frac{1+\alpha}{\alpha} \iint  f_\theta(y|x)^{\alpha}f_{X,\gamma}(x)^{\alpha}d\widehat{G}(x,y) \nonumber\\
	&&= \iint f_\theta(y|x)^{1+\alpha}f_{X,\gamma}(x)^{1+\alpha}dxdy 
	- \frac{1+\alpha}{\alpha} \sum_{i=1}^n W_{in} 
	f_\theta(Z_{(i,n)}|X_{[i,n]})^{\alpha}f_{X,\gamma}(X_{[i,n]})^{\alpha}.~~~~~~
	\label{EQ:obj_fn_MDPDE_SC}
\end{eqnarray} 
For the case $\alpha=0$, the MDPDE of $(\theta,~\gamma)$ is 
to be obtained by minimizing the objective function  
$
\left[\lim\limits_{\alpha \downarrow 0}~H_{n,\alpha}(\theta,~\gamma)\right], 
$
or equivalently 
\begin{eqnarray}
	H_{n, 0}(\theta, \gamma) &=& 
	-\sum_{i=1}^n W_{in} \log\left[f_\theta(Z_{(i,n)}|X_{[i,n]})f_{X,\gamma}(X_{[i,n]})\right].
	\label{EQ:obj_fn_MDPDE_SC0}
\end{eqnarray} 
The estimator obtained by minimizing (\ref{EQ:obj_fn_MDPDE_SC0}) is nothing 
but the maximum likelihood estimator of $(\theta,~\gamma)$ under the present set-up.
Therefore, the proposed MDPDE is indeed a generalization of the MLE.

The estimating equations of the MDPDE of  $(\theta,~\gamma)$ are then given by  
$$
\frac{\partial H_{n,\alpha}(\theta, \gamma)}{\partial\theta} = 0, ~~~~~ 
\frac{\partial H_{n,\alpha}(\theta, \gamma)}{\partial\gamma} = 0. ~~~~~~~\alpha \geq 0.
$$
For $\alpha > 0$, routine differentiation simplifies the estimating equations to yield 
\begin{eqnarray}
	\zeta_\theta - \sum_{i=1}^n W_{in} u_\theta(Z_{(i,n)},X_{[i,n]})f_\theta(Z_{(i,n)}|X_{[i,n]})^{\alpha}f_{X,\gamma}(X_{[i,n]})^{\alpha}
	&=& 0, \label{EQ:est-eqn-1} \\
	\zeta_\gamma - \sum_{i=1}^n W_{in} u_\gamma(X_{[i,n]})f_\theta(Z_{(i,n)}|X_{[i,n]})^{\alpha}f_{X,\gamma}(X_{[i,n]})^{\alpha}
	& =&  0,
	\label{EQ:est-eqn-2}
\end{eqnarray}
where
\begin{eqnarray}
	\zeta_\theta &=& \iint u_\theta(y,x) f_\theta(y|x)^{1+\alpha}f_{X,\gamma}(x)^{1+\alpha}dxdy ,\nonumber\\
	\zeta_\gamma &=& \iint u_\gamma(x) f_\theta(y|x)^{1+\alpha}f_{X,\gamma}(x)^{1+\alpha}dxdy,\nonumber   
\end{eqnarray}
with $u_\theta(y,x) = \frac{\partial \ln f_\theta(y|x)}{\partial\theta}$
and $u_\gamma(x) = \frac{\partial \ln f_\gamma(x)}{\partial\gamma}$ 
being the score functions corresponding to $\theta$ and $\gamma$ respectively.
For $\alpha=0$, the corresponding estimating equation obtained by differentiating $H_{n,0}$ has the simpler 
form given by 
\begin{eqnarray}
	\sum_{i=1}^n W_{in} u_\theta(Z_{(i,n)},X_{[i,n]})&=& 0, \nonumber \\
	\sum_{i=1}^n W_{in} u_\gamma(X_{[i,n]}) & =&  0,\nonumber
\end{eqnarray}
which can also be obtained from Equations (\ref{EQ:est-eqn-1}) and (\ref{EQ:est-eqn-2}) by substituting $\alpha=0$;
note that, at $\alpha=0$, $\zeta_\theta=0$ and $\zeta_\gamma=0$. 
Therefore, Equations (\ref{EQ:est-eqn-1}) and (\ref{EQ:est-eqn-2}) represent the estimating equations for all 
MDPDEs with $\alpha \geq 0$. 

\begin{definition}\label{DEF:MDPDE}
	Consider the above mentioned set-up. The minimum density power divergence estimator of 
	$(\theta,~\gamma)$ based on the observed data $(Z_i, \delta_i, X_i)$, $i=1, \ldots, n$ is defined by the
	simultaneous root of the equations (\ref{EQ:est-eqn-1}) and (\ref{EQ:est-eqn-2}). 
	If there are multiple roots of these equations, the MDPDE will be given by the root 
	which minimizes the objective function (\ref{EQ:obj_fn_MDPDE_SC}) for $\alpha>0$, 
	or (\ref{EQ:obj_fn_MDPDE_SC0}) for $\alpha=0$.
\end{definition}

Clearly, the MDPDE is Fisher consistent by its definition and
the estimating equations (\ref{EQ:est-eqn-1}) and (\ref{EQ:est-eqn-2}) are unbiased at the model.
Further, there will not be any problem of root selection in case of multiple roots, 
which is a general issue in the inferences based on estimating equation. 
This is because we have a proper objective function in the case of MDPDE.

Note that the MDPDE estimating equations (\ref{EQ:est-eqn-1}) and (\ref{EQ:est-eqn-2}) 
can also be written as
\begin{eqnarray}
	\iint  \psi(y, x; \theta, \gamma) d\widehat{G}(x,y) 
	=\sum_{i=1}^n W_{in} \psi(Z_{(i,n)},X_{[i,n]}; \theta, \gamma)  = 0, 
	\label{EQ:est-eqn-psi} 
\end{eqnarray}
where $\psi(y, x; \theta, \gamma)  = (\psi_{1,\alpha}(y, x; \theta, \gamma) , ~ \psi_{2,\alpha}(y, x; \theta, \gamma))^T$ with 
\begin{eqnarray}
	\left.\begin{array}{rcl}
		\psi_{1,\alpha}(Y, X; \theta, \gamma)  &=& \zeta_\theta  
		-  u_\theta(Y, X) f_\theta(Y|X)^{\alpha}f_{X,\gamma}(X)^{\alpha},  \\ 
		\psi_{2,\alpha}(Y, X; \theta, \gamma) &=& \zeta_\gamma
		- u_\gamma(X)f_\theta(Y|X)^{\alpha}f_{X,\gamma}(X)^{\alpha}. 
	\end{array}\right\}\label{EQ:psi-2}
\end{eqnarray}
This particular estimating equation (\ref{EQ:est-eqn-psi}) is similar to 
that of the M-estimator for i.i.d.~non-censored data with covariates 
(in fact they become the same if $\widehat{G}$ is the empirical distribution function).
So, extending the concept of Wang (1999),  we can define the general M-estimator of  
$\theta$ and $\gamma$ based on any general  function 
\begin{eqnarray}
	\psi(y, x; \theta, \gamma)~~:~~\mathbb{R}\times \mathbb{R}^p \times \mathbb{R}^q \times \mathbb{R}^r 
	\mapsto \mathbb{R}^q \times \mathbb{R}^r \label{EQ:psi-prop0}
\end{eqnarray}
as the solution of the estimating equation (\ref{EQ:est-eqn-psi}).
However, to make it an unbiased estimating equation, we consider only the $\psi$-functions for which 
\begin{eqnarray}
	\iint  \psi(y, x; \theta, \gamma) d{G}(x,y) =0.
	\label{EQ:psi-prop1}
\end{eqnarray}

\begin{definition}\label{DEF:M-est}
	Consider the above mentioned parametric set-up for censored data with stochastic covariates.
	Also, consider a general $\psi$-function as in (\ref{EQ:psi-prop0}) satisfying the condition (\ref{EQ:psi-prop1}).
	An M-estimator of $(\theta,~\gamma)$  corresponding to this general $\psi$-function based on the observed data 
	$(Z_i, \delta_i, X_i)$, $i=1, \ldots, n$, is defined as  the root of the 
	estimating equation (\ref{EQ:est-eqn-psi}). 
\end{definition}

Note that any general M-estimator is also Fisher consistent and is based on an unbiased estimating equation, 
by definition. But, in general, they may suffer from the problem of multiple roots and 
need a proper numerical techniques (like bootstrapping) to get a well-defined M-estimator. 
So, in this paper, we restrict our attention mainly to examine the performances of the robust MDPDE 
corresponding to the particular $\psi$-function defined in (\ref{EQ:psi-2}).
However, we derive theoretical asymptotic properties of the general M-estimators in Section \ref{SEC:M-estimators} 
and deduce the properties of MDPDEs from those general results in Section \ref{SEC:MDPDE_prop}. 

%
%


\subsection{Application (I): Fully Parametric version of Linear Regression model }
\label{SEC:regression}

We first consider the simplest problem of linear regression with censored responses and stochastic covariables. 
Precisely, we assume the linear regression model (LRM)
\begin{equation}
	Y_i = X_i^T\theta + \epsilon_i, ~~~~ i=1, \ldots, n,
\end{equation}
where $Y_i$ is the censored response (generally the lifetime), $X_i$ is a $p$-variate 
stochastic auxiliary variable associated with the response, $\theta$ is the vector of unknown
regression coefficients and $\epsilon_i$ is the error in specified linear model. 
We assume that the error  $\epsilon_i$s are independent and identically distributed 
with distribution function $F_{e}$ 
and $X_i$s are independent of the errors having distribution function $F_{X,\gamma}$.
Generally we can assume both symmetric error distributions like the normal as well as 
asymmetric error distributions like the exponential;
however the second group is used in most reliability applications. 
Then, the conditional distribution of the response variable $Y_i$ given $X_i$ is 
$F_\theta(y|X_i) = F_e(y- X_i^T\theta)$.
Now, we consider the incomplete censored observations $(Z_i, \delta_i)$ as defined in Section 1
and use them to estimate $(\theta, \gamma)$ robustly and efficiently.

This inference problem clearly belongs to the general set-up considered in the previous subsection and 
frequently arises in reliability studies and other applied researches. 
We can obtain a robust solution to this problem through the proposed MDPDEs, 
obtained by just solving the estimating equations (\ref{EQ:est-eqn-1}) and (\ref{EQ:est-eqn-2}).
Here we present the detail working for one particular example of model families 
$\mathcal{F_X}$ and $\mathcal{F_0}$. The case of other model families can also be tackled similarly.

Suppose the response variable is exponentially distributed with mean depending 
on the covariates as $E[Y|X] = X^T\theta$; then 
$\mathcal{F_X} = \left\{ Exp(X^T\theta) : \theta \in \mathbb{R}^p \right\}$ 
where $Exp(\tau)$ represents the exponential distribution with mean $\tau$. 
Also, for simplicity, let us assume that the auxiliary variables are independent to each other and 
normally distributed so that $\mathcal{F_0}=\left\{N_p(\gamma, I_p) : \gamma \in \mathbb{R}^p\right\}$. 
In this case, the objective function $H_{n,\alpha}(\theta, \gamma)$ of the MDPDE has a simpler form given by 
\begin{eqnarray}
	H_{n,\alpha}(\theta, \gamma) 
	&=& \frac{(1+\alpha)^{-3/2}}{(2\pi)^{\alpha/2}} \psi^{(0)}(\theta,\gamma)
	-\frac{1+\alpha}{(2\pi)^{\alpha/2}\alpha} \sum_{i=1}^n W_{in}\frac{e^{-\alpha\psi_i(\theta,\gamma)}}{(X_{[i,n]}^T\theta)^{\alpha}} 
	~~~\alpha>0, ~~~~~
	\label{EQ:obj_fn_MDPDE_Reg}\\
	\mbox{and   }~~~
	H_{n, 0}(\theta, \gamma) &=&  \sum_{i=1}^n W_{in}
	\left[\psi_i(\theta,\gamma) + \log(X_{[i,n]}^T\theta) + \frac{1}{2}\log(2\pi)\right],
	\label{EQ:obj_fn_MDPDE_Reg0}
\end{eqnarray}
where
$\psi^{(0)}(\theta,\gamma)  = \int (x^T\theta)^{-\alpha} N_p(x, \gamma, \frac{1}{1+\alpha}I_p) dx 
$
and 
$\psi_i(\theta,\gamma)= \frac{Z_{(i,n)}}{X_{[i,n]}^T\theta} 
+ \frac{1}{2}(X_{[i,n]} - \gamma)^T(X_{[i,n]} - \gamma)$. 
Note that the integral $\psi^{(0)}(\theta,\gamma)$  is just the expectation of 
a simple function of multivariate normal random variable;
so it can be computed  quite easily using standard numerical integration techniques. 
Therefore, we can simply minimize the above objective functions by any numerical algorithm to obtain the MDPDE at any $\alpha \geq 0$. 

Alternatively, we can obtain the MDPDE by solving the estimating equations 
(\ref{EQ:est-eqn-1}) and (\ref{EQ:est-eqn-2}), which can also be simplified in this particular 
situation as
\begin{eqnarray}
	\frac{\alpha}{(1+\alpha)^{5/2}} \bar{\psi}^{(0)}(\theta,\gamma) &=&
	\sum_{i=1}^n W_{in} e^{-\alpha\psi_i(\theta,\gamma)}
	\left(\frac{1}{(X_{[i,n]}^T\theta)^{(1+\alpha)}} - \frac{Z_{(i,n)}}{{(X_{[i,n]}^T\theta)^{2}}}\right) X_{[i,n]},
	\nonumber \\
	\frac{\alpha}{(1+\alpha)^{5/2}} \left[\bar{\psi}^{(0)}(\theta,\gamma) - \gamma {\psi}^{(0)}(\theta,\gamma)\right]
	&=& \sum_{i=1}^n W_{in} e^{-\alpha\psi_i(\theta,\gamma)}
	\frac{(X_{[i,n]} - \gamma)}{(X_{[i,n]}^T\theta)^{\alpha}},
	\nonumber
\end{eqnarray}
where
$
\bar{\psi}^{(0)}(\theta,\gamma)  = \int (x^T\theta)^{-\alpha} x N_p(x, \gamma, \frac{1}{1+\alpha}I_p) dx. 
$
The case of $\alpha=0$ (MLE) can be simplified further 
where the estimator of $\gamma$ becomes independent of the parameter $\theta$. 
To see this, we simplify the above estimating equations at $\alpha=0$ as
\begin{eqnarray}
	\sum_{i=1}^n W_{in} \frac{\left(Z_{(i,n)} - X_{[i,n]}^T\theta\right)}{{(X_{[i,n]}^T\theta)^{2}}} X_{[i,n]} = 0,
	~~~~~~~\sum_{i=1}^n W_{in}\left(X_{[i,n]} - \gamma\right) =0. \nonumber
\end{eqnarray}
Solving the second equation, we get that 
$
\widehat{\gamma} = \frac{\sum_{i=1}^n W_{in}X_{[i,n]}}{\sum_{i=1}^n W_{in}},
$
which clearly does not depend on $\theta$, as is expected from the theory of maximum likelihood inference.

The semi-parametric version of this model has been considered in Zhou (2010), where no assumption has been made about 
the distribution of the i.i.d. sequences $\{X_i\}$ and $\{\epsilon_i\}$.
In that paper, an M-estimator of $\theta$ has been proposed by solving the estimating equation 
\begin{equation}
	\sum_{i=1}^n W_{in} \psi_0\left(Z_{(i,n)} - X_{[i,n]}^T\theta\right) = 0,
	\label{EQ:EstEqn_Zhou}
\end{equation}
for suitable choices of $\psi_0$. Although this M-estimator is different from the proposed MDPDE,
it in fact belongs to our general class of M-estimators (Definition \ref{DEF:M-est}) as we will show in detail in section \ref{SEC:SemiP_Ex}.
Further, we will show in Section \ref{SEC:IF} that our MDPDE with $\alpha>0$ under fully parametric set-up 
is robust with respect to both the outliers in response and the leverage points in the explanatory variables;
but the M-estimator of Zhou (2010) is not robust with respect to the leverage points.

\subsection{Application (II): Fully Parametric Exponential Regression model for Medical Sciences}
\label{SEC:Exp_regression}

The simple linear regression model considered in the previous subsection is  
the most popular inference problem under the set-up considered in this paper.
However, although it is simple and potentially applicable in several real life problems, 
for the purpose of serving the typical applications in the medical science this simple linear model is rarely used. 
The reason is that in almost all medical applications the support of the distribution of 
the censoring times is shorter than the support of the lifetimes (we cannot follow the patients until they die). 
For this reason the linear model is usually not identifiable and 
in order to become applicable the setup of the proposed method
needs to introduce a truncation time, say $\tau$, such that the probability to be uncensored by time $\tau$ 
is strictly greater than zero for all $x$.
We can suitably extend the proposed MDPDE to cover these assumptions through some more routine calculations,
which we leave for the readers.

In this section we present an alternative multiplicative model with exponential error for the applications in medical sciences. 
This particular model, known as the exponential regression model, is widely used 
and most popular in the medical sciences and related applications.
More precisely, let us assume the multiplicative regression model for the survival times (responses) 
$Y_i$, $i=1, \ldots,n$ as 
\begin{equation}
	Y_i = e^{X_i^T\theta} \times \epsilon_i, ~~~~ i=1, \ldots, n,
\end{equation}
where $X_i$ is a $p$-variate stochastic auxiliary variable associated with the response, 
$\theta$ is the vector of unknown regression coefficients and $\epsilon_i$ is the error in the specified linear model. 
Such a multiplicative model ensures the positivity of the the response variables, 
which are generally life-time in most applications.
In the exponential regression model, we assume that the error variable $\epsilon$ 
is exponentially distributed with mean $1$. 
Then the conditional distribution of the response variable $Y$ given the covariate $X$ 
is also exponential with mean $E[Y|X] = e^{X^T\theta}$; 
so considering the notations of Section \ref{SEC:intro}, 
$\mathcal{F_X} = \left\{ Exp(e^{X^T\theta}) : \theta \in \mathbb{R}^p \right\}$. 
Also, the  $X_i$s are independent of the errors having distribution function $F_{X,\gamma}$.
Our objective is to estimate $(\theta, \gamma)$ robustly and efficiently based on the 
incomplete (censored) observations $(Z_i, \delta_i, X_i)$ as defined in Section \ref{SEC:intro}.

Again this inference problem belongs to the general set-up of Section \ref{SEC:general_MDPDE}
so that the proposed MDPDEs provide a robust solution to it. 
In the case of independent and normally distributed covariates with 
$\mathcal{F_0}=\left\{N_p(\gamma, I_p) : \gamma \in \mathbb{R}^p\right\}$, 
we can simplify the objective function $H_{n,\alpha}(\theta, \gamma)$, 
to be  minimized in order to obtain the MDPDE, 
as
\begin{eqnarray}
	H_{n,\alpha}(\theta, \gamma) 
	&=& \frac{e^{\alpha(\gamma^T\theta) + \frac{\alpha^2}{2(1+\alpha)}(\theta^T\theta)}}{(1+\alpha)^{3/2}(2\pi)^{\alpha/2}} 
	-\frac{(1+\alpha)}{(2\pi)^{\alpha/2}\alpha} \sum_{i=1}^n W_{in} e^{-\alpha\Gamma_i(\theta,\gamma)}
	~~~\alpha>0, ~~~~~
	\label{EQ:obj_fn_MDPDE_ExpReg}\\
	\mbox{and   }~~~~
	H_{n, 0}(\theta, \gamma) &=&  \sum_{i=1}^n W_{in}
	\left[ \Gamma_i(\theta,\gamma) + \frac{1}{2}\log(2\pi)\right],
	\label{EQ:obj_fn_MDPDE_ExpReg0}
\end{eqnarray}
where $\Gamma_i(\theta,\gamma) = {Z_{(i,n)}}e^{X_{[i,n]}^T\theta} + (X_{[i,n]}^T\theta) 
+ \frac{1}{2}(X_{[i,n]} - \gamma)^T(X_{[i,n]} - \gamma).
$
This objective function can be easily minimized using any standard numerical techniques for any $\alpha \geq 0$. 

The estimating equations  (\ref{EQ:est-eqn-1}) and (\ref{EQ:est-eqn-2}) can also be simplified 
for this particular case of ERM with normal covariates. For $\alpha >0$, they have the form
\begin{eqnarray}
	\sum_{i=1}^n W_{in}  e^{-\alpha\Gamma_i(\theta,\gamma)}
	\left(Z_{(i,n)} e^{(X_{[i,n]}^T\theta)} +1\right) X_{[i,n]} &=& \frac{\alpha}{(1+\alpha)^{5/2}}  
	\left[\gamma + \frac{\alpha\theta}{(1+\alpha)}\right]
	e^{\alpha(\gamma^T\theta) + \frac{\alpha^2}{2(1+\alpha)}(\theta^T\theta)},
	\nonumber\\
	\sum_{i=1}^n W_{in}  e^{-\alpha\Gamma_i(\theta,\gamma)}
	\left(X_{[i,n]}- \gamma\right) &=& \frac{\alpha^2\theta}{(1+\alpha)^{5/2}}  
	e^{\alpha(\gamma^T\theta) + \frac{\alpha^2}{2(1+\alpha)}(\theta^T\theta)}.
	\nonumber
\end{eqnarray}
At $\alpha=0$ (MLE), these estimating equations further simplifies to
\begin{eqnarray}
	\sum_{i=1}^n W_{in} \left(Z_{(i,n)} e^{(X_{[i,n]}^T\theta)} +1\right) X_{[i,n]}  = 0,
	~~~~~~ \sum_{i=1}^n W_{in}\left(X_{[i,n]} - \gamma\right) =0, \nonumber
\end{eqnarray}
which again produce the same estimator of $\gamma$ as in the case of the LRM and independent of the 
parameter $\theta$.

\subsection{Application (III): Fully Parametric Version of Accelerated Failure Time (AFT) Models}
\label{SEC:AFT}

The exponential regression model considered in the previous section can be linearized by taking natural logarithm
of the response time:
\begin{equation}
	\log(Y_i) = X_i^T\theta + \epsilon_i^*,
\end{equation}
where $\epsilon_i^* =\log(\epsilon)$ follows the standard extreme value distribution.
This model can be generalized by considering some alternative distribution for $\epsilon_i^*$, but with mean 0.
When $\epsilon_i^*$ follows an extreme value distribution with mean $0$ and scale parameter $\sigma$,
then $Y_i$ follows a Weibull distribution and the resulting regression model is known as the Weibull regression model (WRM).
Other common distributions for $\epsilon_i^*$ are logistic (survival time has log-logistic distribution),
normal (survival times are log-normal) etc.
In such models, the covariate has a multiplicative effect of the response life-time and hence 
they are generally known as the accelerated failure time (AFT) model.

Consider a general location scale model family $\left\{f_{\mu,\sigma}(x)=f((x-\mu)/\sigma)\right\}$ 
for some known function $f$
and let $\epsilon_i^*$ have density $f_{0,\sigma}$. Then, given the covariate $X$,  $Y_i$ has density 
$$
\frac{1}{y}f\left(\frac{log(y)-(X^T\beta)}{\sigma}\right).
$$ 
This is the general form of the parametric AFT regression model;
taking $f$ as standard extreme value distribution it simplifies to the Weibull regression model and so on.
Suppose the distribution of the covariates $X_i$ is modeled by the family $F_{X,\gamma}$ as in the earlier cases
so that our target becomes the estimation of $\theta=(\beta, \sigma)$ and $\gamma=\mu$ based on the censored observations 
$(Z_i, \delta_i, X_i)$ as defined in Section \ref{SEC:intro}.

Now we can again minimize the objective function (\ref{EQ:obj_fn_MDPDE_SC}) 
with respect to the parameters to obtain their robust MDPDE.
The exact form of this objective function can be obtained easily for any particular choice of $f$
and any standard numerical algorithm provide us with the solution of this optimization problem.

Further, as we will see through the numerical illustrations in Section \ref{SEC:Numerical},
the proposed  MDPDE provides highly robust solutions in presence of outliers in data with little loss in efficiency
at small positive $\alpha$. In the case of the accelerated failure time models, the robustness of the MDPDE is directly 
comparable with the alternative proposal of Locatelli et al.~(2011). 
Note that, contrary to the Locatelli et al.~(2011) approach, our proposed estimator can estimate the parameter ($\gamma $)
in the distribution of covariates simultaneously with $\theta=(\beta, \sigma)$; 
although someone might not see it to be a big advantage as the parameter $\gamma$ can be estimated separately in many cases. 
However, the major advantage of our proposal is its generality and computational simplicity
for any kind of parametric model with censored survival data with stochastic covariates.
We will see in Section \ref{SEC:SemiP_Ex} that, even if we ignore the estimation of the parameters of the marginal
distribution of covariates, our proposal contains
the existing proposals of Locatelli et al.~(2011), a similar proposal of Zhou (2010) under the same set-up
and also their extension in Wang et al.~(2015).
This vast generality is the major strength of our proposal over the existing literature; 
see Section \ref{SEC:SemiP_Ex} for more detailed and general discussions.

\section{Asymptotic Properties}
\label{SEC:M-estimators}

Consider the models and  set-up described in Section \ref{SEC:general_MDPDE}. 
First, we derive the asymptotic properties of the general M-estimator $(\hat\theta_n, \hat{\gamma}_n)$
of $(\theta, \gamma)$ as defined in Definition \ref{DEF:M-est} based on a (random) censored sample of size $n$. 
Let us assume that the true distributions belong to the corresponding model families
with $(\theta_0, \gamma_0)$ being the true parameter value. 
Define 
$$
\lambda_n(\theta, \gamma) = \iint  \psi(y, x; \theta, \gamma) d\widehat{G}(x,y) 
=\sum_{i=1}^n W_{in} \psi(Z_{(i,n)},X_{[i,n]}; \theta, \gamma),
$$
and 
$$
\lambda_G(\theta, \gamma) = \iint  \psi(y, x; \theta, \gamma) d{G}(x,y). 
$$

Then $\lambda_n(\theta, \gamma)$ is the empirical version of $\lambda_G(\theta, \gamma)$; also
by definition $\lambda_n(\hat\theta_n, \hat{\gamma}_n)=0$ and $\lambda_G(\theta_0, \gamma_0) = 0$.
In order to prove the asymptotic consistency and normality of the general M-estimator 
$(\hat\theta_n, \hat{\gamma}_n)$, we use the results of Section \ref{SEC:Background}.
So, we will assume that the assumptions (A1) to (A4) holds true with $\phi$ replaced by $\psi$. 
Further, let us assume

\begin{itemize}
	\item[(A5)] Either $\tau_{G_Y} < \tau_{G_C}$ or $\tau_{G_C}=\infty$ so that $\widetilde{G}$ and $G$ 
	coincides.
	
	\item[(A6)] The  variance matrix $\Sigma_\psi$, as defined in (\ref{EQ:CLT_var})
	with $\phi$ replaced by $\psi$, exists finitely (with all entries finite).
\end{itemize}
Then Proposition \ref{PROP:Consistency} and  Proposition \ref{PROP:CLT}
give the strong consistency and asymptotic normality of $\lambda_n(\theta, \gamma)$, 
which are summarized in the following lemma.

\begin{lemma}
	Consider the above set-up with an integrable function $\psi$ and 
	let $(\theta_0, \gamma_0)$ be the true parameter value. Then,
	\begin{itemize}
		\item[(i)] Under Assumptions (A1), (A2) and (A5) 
		we have, with probability one, 
		$$
		\lim\limits_{n\rightarrow\infty} \lambda_n(\theta_0, \gamma_0)  = \lambda_G(\theta_0, \gamma_0)
		$$
		
		\item[(ii)] Under Assumptions (A2) to (A6), the asymptotic distribution of 
		$$
		\sqrt{n} \left[\lambda_n(\theta, \gamma) - \lambda_G(\theta, \gamma)\right] 
		= \sqrt{n} \iint  \psi(y, x; \theta, \gamma) d[\widehat{G}-G](x,y)
		$$
		is normal with mean $0$ and variance matrix $\Sigma_\psi$.
	\end{itemize}
	\label{LEM:asym_lambda}
\end{lemma}

\subsection{Strong Consistency and Asymptotic Normality of General M-estimators}
\label{SEC:SC_AN_M-est}

Wang (1999) has proved the strong consistency and asymptotic normality results for the M-estimators 
based on only censored variable with no covariables. This section extend the theory to the case 
where covariables are present along with the censored response. 
The extensions are in the line of the corresponding results with no censoring (see Huber, 1981; Serfling, 1980).
Further note that whenever $\gamma$ is known the asymptotic properties of $\hat\theta_n$ 
follows from just a routine application of the results derived in Wang (1999); 
so here we assume $\gamma$ to be known and derive the joint distribution of M-estimator of $\theta$ and $\gamma$.
These results provide a general (asymptotic) theoretical framework 
to study the properties of a wide class of estimators of $(\theta, \gamma)$ depending on the estimating equations.

Denote the $j$-th component of the function $\psi$ by $\psi_j$ for $j=1, \ldots, q+r$.
Also, consider following (stronger) conditions on the nature of the function $\psi$.
\begin{itemize}
	\item[(A7)] $\psi(y, x; \theta, \gamma)$ is continuous in $(\theta, \gamma)$ and also bounded.
	
	\item[(A8)] The population estimating equation $\lambda_G(\theta, \gamma)=0$ 
	has an unique root given by $(\theta_0, \gamma_0)$.
	
	\item[(A9)] There exists a compact set $C$ in $\mathbb{R}^q \times \mathbb{R}^r$ satisfying
	$$
	\inf_{(\theta, \gamma) \notin C} \left|\iint \psi_j(y, x; \theta, \gamma) dG(x, y)\right| > 0, 
	~~~~~~~~ j=1, \ldots, q+r.
	$$
\end{itemize}
Now let us start with the strong consistency of the M-estimator $(\hat\theta_n, \hat{\gamma}_n)$
by an extension of Theorem 3 of Wang (1999, page 307) under above conditions. 
The proof follows in the same line of Wang (1999) by replacing the corresponding SLLN,
given in Proposition 1 of Wang (1999), by the part (i) of Lemma \ref{LEM:asym_lambda} in the present context;
hence it is omitted for simplicity of presentation.

\begin{theorem}
	Consider the above set-up with Assumptions (A1), (A2), (A5), (A7) and (A8).
	Then we have the following results.
	\begin{itemize}
		\item[(i)] There exists a sequence of M-estimators $\{(\hat{\theta}_n, \hat{\gamma}_n)\}$
		satisfying the empirical estimating equation $\lambda_n(\theta, \gamma)=0$
		that converges with probability one to $(\theta_0, \gamma_0)$. 
		
		\item[(ii)] Further if (A9) also holds true, then any sequence of M-estimators 
		$\{(\hat{\theta}_n, \hat{\gamma}_n)\}$ satisfying $\lambda_n(\theta, \gamma)=0$
		converges with probability one to $(\theta_0, \gamma_0)$.
	\end{itemize}
	\label{THM:SLLN_M-est}
\end{theorem}

Note that the first part (i) of Theorem \ref{THM:SLLN_M-est} 
is just a multivariate extension of Lemma B of Serfling (1980, page 249)
from the complete data case to the present case of censored data with covariates.
Further, the additional condition (A9) in part (ii) makes any sequence of M-estimators 
satisfying the estimating equation (\ref{EQ:est-eqn-psi}) to eventually fall 
in a compact neighborhood of $(\theta_0, \gamma_0)$. 
This result, even with the stronger conditions, becomes really helpful when the 
empirical estimating equation $\lambda_n(\theta, \gamma)=0$ has multiple roots 
and one could obtain different M-estimator sequences by applying 
different numerical equation solving techniques. 
Part (ii) of Theorem \ref{THM:SLLN_M-est} ensures that all theses sequences 
of M-estimators will be strongly consistent for the unique root 
$(\theta_0, \gamma_0)$ of the equation $\lambda_G(\theta, \gamma)=0$.

Next we turn our attention to the asymptotic normality of M-estimators. 
In this regard, we will first present a useful lemma in terms 
of any real valued function $g(y, x, \theta_0, \gamma_0)$. 
This is again a suitable extension of Lemma 1 of Wang (1999, page 307)
to the present set-up and the proof follows similarly by replacing the corresponding SLLN 
(Proposition 1 of Wang) by Part (i) of Lemma \ref{LEM:asym_lambda}. 
Assume the following condition about the function   $g(y, x, \theta_0, \gamma_0)$.

\begin{itemize}
	\item[(A10)] For a  real valued function $g(y, x, \theta_0, \gamma_0)$, at least one of the following holds:
	\begin{itemize}
		\item[(i)] $g(y, x, \theta, \gamma)$ is continuous at $(\theta_0, \gamma_0)$ uniformly in $(y,x)$.
		
		\item[(ii)] As $\delta \rightarrow 0$, 
		$$\iint \sup_{\{(\theta, \gamma) : ||(\theta, \gamma) - (\theta_0, \gamma_0)|| \leq \delta\}} 
		\left| g(y, x, \theta, \gamma) - g(y, x, \theta_0, \gamma_0)\right|dG(x,y) = h_\delta \rightarrow 0.
		$$
		(Here $||\cdot||$ denotes the Euclidean norm).
		
		\item[(iii)] $g$ is continuous in $(y,x)$ for for any fixed $(\theta, \gamma)$ in a neighborhood of $(\theta_0, \gamma_0)$,
		and 
		$$
		\lim\limits_{(\theta, \gamma) \rightarrow (\theta_0, \gamma_0)}
		\left|\left| g(y, x, \theta, \gamma) - g(y, x, \theta_0, \gamma_0)\right|\right|_v  = 0.
		$$
		(Here $||\cdot||_v$ denotes the total variation norm).
		
		\item[(iv)] $\iint g(y, x, \theta, \gamma) dG(x,y) $ is continuous at $(\theta, \gamma) = (\theta_0, \gamma_0)$,
		and $g$ is continuous in $(y,x)$ for $(\theta, \gamma)$ in a neighborhood of $(\theta_0, \gamma_0)$,
		and 
		$$
		\lim\limits_{(\theta, \gamma) \rightarrow (\theta_0, \gamma_0)}
		\left|\left| g(y, x, \theta, \gamma) - g(y, x, \theta_0, \gamma_0)\right|\right|_v < \infty.
		$$
		
		\item[(v)]  $\iint g(y, x, \theta, \gamma) dG(x,y) $ is continuous at $(\theta, \gamma) = (\theta_0, \gamma_0)$,
		and 
		$$
		\iint g(y,x,\theta,\gamma)d\hat{G}(x,y)\mathop{\rightarrow}^\mathcal{P}\iint g(y,x,\theta,\gamma)dG(x,y)<\infty,
		$$
		uniformly for $(\theta, \gamma)$ in a neighborhood of $(\theta_0, \gamma_0)$. 
	\end{itemize}
\end{itemize}

\begin{lemma}
	Suppose $g(y, x, \theta_0, \gamma_0)$ is a real valued function with 
	$\iint g(y, x, \theta_0, \gamma_0) dG(x,y) < \infty$. 
	Assume that the conditions (A1), (A2) and (A10) hold for $g$. Then, for any sequence 
	$(\hat{\theta}_n, \hat{\gamma}_n) \displaystyle\mathop{\rightarrow}^\mathcal{P} (\theta_0, \gamma_0)$, 
	we have 
	$$
	\iint g(y,x,\hat{\theta}_n, \hat{\gamma}_n)d\hat{G}(x,y) \mathop{\rightarrow}^\mathcal{P}
	\iint g(y,x,\theta_0,\gamma_0)dG(x,y).
	$$
	\label{LEM:asym_lem}
\end{lemma}

\begin{theorem}
	Consider the above set-up and assume that $\psi$ is differentiable with respect to 
	$(\theta, \gamma)$ in a neighborhood of $(\theta_0, \gamma_0)$ and the matrix
	\begin{eqnarray}
		\Lambda_G(\theta_0, \gamma_0) = \iint 
		\left.\frac{\partial}{\partial(\theta, \gamma)}\psi(y,x,\theta,\gamma)\right|_{(\theta, \gamma)=(\theta_0, \gamma_0)}
		dG(x,y),
	\end{eqnarray}
	exists finitely and is non-singular.  
	Further assume that the assumptions of Lemma \ref{LEM:asym_lem} hold for 
	$g(y,x,\theta,\gamma)=\Lambda_G^{ij}(\theta_0,\gamma_0)$, the $(i,j)$-th element of $\Lambda_G(\theta_0,\gamma_0)$, 
	with $i, j=1, \ldots, q+r$. 
	Then, under Assumptions (A2) to (A6) we have, for any sequence of M-estimators 
	$\{(\hat{\theta}_n, \hat{\gamma}_n)\}$ satisfying $\lambda_n(\theta, \gamma)=0$
	that converges in probability to $(\theta_0, \gamma_0)$, 
	$$
	\sqrt{n} \left[(\hat{\theta}_n, \hat{\gamma}_n) - (\theta_0, \gamma_0) \right] \mathop{\rightarrow}^\mathcal{D} 
	N\bigg(0, \Lambda_G(\theta_0, \gamma_0)^{-1}\Sigma_\psi(G)\Lambda_G(\theta_0, \gamma_0)^{-1}\bigg).
	$$
	\label{THM:AsymN_M-est}
\end{theorem}
\noindent
\textbf{Proof: }
Since $\psi$ is differentiable in $(\theta, \gamma)$ , so  is the function $\lambda_n(\theta, \gamma)$.
So an application of multivariate mean value theorem yields 
$$
\lambda_n(\hat{\theta}_n, \hat{\gamma}_n) - \lambda_n(\theta_0, \gamma_0) 
= \Lambda_{\hat{G}}(\zeta_{1n}, \zeta_{2n})\left[(\hat{\theta}_n, \hat{\gamma}_n) - (\theta_0, \gamma_0)\right],
$$ 
with 
$||(\zeta_{1n}, \zeta_{2n}) - (\theta_0, \gamma_0)||< ||(\hat{\theta}_n, \hat{\gamma}_n) - (\theta_0, \gamma_0)||$.
Further, by definition, $\lambda_n(\hat{\theta}_n, \hat{\gamma}_n)=0$ and $\lambda_G(\theta_0, \gamma_0) =0$.
Hence we get, 
$$
(\hat{\theta}_n, \hat{\gamma}_n) - (\theta_0, \gamma_0) 
= - \left[\Lambda_{\hat{G}}(\zeta_{1n}, \zeta_{2n})\right]^{-1} 
\left[\iint  \psi(y, x; \theta, \gamma) d[\widehat{G}-G](x,y)\right].
$$ 
However, it follows from Lemma \ref{LEM:asym_lem} that  
each term of $\Lambda_{\hat{G}}(\zeta_{1n}, \zeta_{2n})$  convergence in probability 
to the corresponding term of $\Lambda_{G}(\theta_0, \gamma_0)$. 
Then, an application of Slutsky's theorem and Part (ii) of Lemma \ref{LEM:asym_lambda} 
completes the proof of the theorem.
\hfill{$\square$}
\bigskip

It is to be noted that the asymptotic normality of the M-estimators require more conditions than 
that required for its strong consistency in terms of differentiability properties of the $\psi$ function,
but it avoid the strong assumptions (A7) -- (A9) used in Theorem \ref{THM:SLLN_M-est}. 
In fact, to obtain the asymptotic distributional convergence of any sequence of M-estimators in this case, 
it is just enough to ensure their convergence to the true parameter value in probability.
All the related conditions used here are in the same spirit with that used in Wang (1999)
have been no covariables are present and were discussed extensively in that paper.

Finally, note that the estimating equation of any general M-estimator can be solved through an appropriate 
numerical technique but the complexity in terms of the iterative procedure increases extensively for 
a complicated non-linear $\psi$-function. However, one can show that, for the Newton-Raphson algorithm, 
if we start the iterations with some $\sqrt{n}$-consistent estimator of $(\theta, \gamma)$ 
then the estimator obtained by just one iteration, known as the one-step M-estimator, 
will have the same asymptotic distribution as the fully iterated M-estimator even 
in case of censored data with covariables as considered here. 
This is a well-known property of the M-estimator in case of complete data.
The following theorem present this precisely for our case; 
the proof follows by an argument similar to that of Theorem 6 of Wang (1999) replacing 
Proposition 1 and 2 of that paper by Part (i) and Part (ii) of Lemma \ref{LEM:asym_lambda}  respectively.

\begin{theorem}
	Suppose the conditions of Theorem \ref{THM:AsymN_M-est} hold true and 
	let $(\widetilde{\theta}_n, \widetilde{\gamma}_n)$ is any $\sqrt{n}$-consistent estimate 
	of the true parameter value $(\theta_0, \gamma_0)$. Then, the one-step M-estimator 
	$(\theta_n^{(1)}, {\gamma}_n^{(1)})$, defined as
	\begin{eqnarray}
		(\theta_n^{(1)}, {\gamma}_n^{(1)}) = (\widetilde{\theta}_n, \widetilde{\gamma}_n)
		- \left[\Lambda_{\hat{G}}(\widetilde{\theta}_n, \widetilde{\gamma}_n)\right]^{-1}
		\lambda_n(\widetilde{\theta}_n, \widetilde{\gamma}_n),
		\label{DFN:One-step_M-est}
	\end{eqnarray}
	has the same distribution as that of the M-estimator $(\hat{\theta}_n, \hat{\gamma}_n)$ 
	derived in Theorem \ref{THM:AsymN_M-est}. 
	\label{THM:One-step_M-est}
\end{theorem}


\subsection{Properties of the MDPDE}\label{SEC:MDPDE_prop}

Note that, the MDPDE is a particular M-estimator with the $\psi$-function given by (\ref{EQ:psi-2})
and so all the results derived in the previous subsection for general  M-estimators also 
hold true for the MDPDEs. In particular MDPDEs are strongly consistent and asymptotically normal
under the assumptions considered in Theorems \ref{THM:SLLN_M-est} and \ref{THM:AsymN_M-est}.
However, in this particular case of MDPDEs, we can closely investigate the required assumptions 
for the particular form of the $\psi$-function.

Note that assumptions (A1), (A2) and (A5) are related to the censoring scheme under consideration
and others are about the special structure of the $\psi$-function.
Further, in this particular case of MDPDE, the $\psi$-function depends on the model density and 
its score function. So, conditions (A3), (A4) and (A6) can easily be shown to hold 
for most statistical models by using the existence of finite and continuous second order moments 
of the score functions with respect to the true distribution $G$.
Similar differentiability conditions on the model and score functions further ensure 
the assumptions of Lemma \ref{LEM:asym_lem}. 
So, the asymptotic normality of the MDPDEs follows from Theorem \ref{THM:AsymN_M-est}
for most models provided we can prove its consistency.
However, assumptions (A6)--(A9), required to prove the strong consistency in Theorem \ref{THM:SLLN_M-est},
are rather difficult one and may not always hold for the assumed model.

Noting that, the asymptotic normality of MDPDEs, as obtained in Theorem \ref{THM:AsymN_M-est}, 
does not require its strong consistency (only convergence in probability is enough),
we now present an alternative approach to prove the (weak) consistency for the particular case of 
MDPDEs under some simpler conditions. This approach is essentially due ot Lehmann (1983), and has been used by Basu et al.~(1998) to 
prove the asymptotic properties of the  MDPDEs under i.i.d.~complete data and extended by 
many researchers later in the context of different  inference problems. 
Here, we extend their approach further for the present case of censored data with covariates. 
Let us also relax the assumption that the true distribution $G$ belongs to the model family
in the sense of assumption (D1) below.
Define 
\begin{eqnarray}
	V(Y, X; \theta, \gamma)  &=& \iint f_\theta(y|x)^{1+\alpha}f_{X,\gamma}(x)^{1+\alpha}dxdy  
	-  \frac{1+\alpha}{\alpha} f_\theta(Y|X)^{\alpha}f_{X,\gamma}(X)^{\alpha}, \nonumber 
	\label{EQ:V_theta} 
\end{eqnarray}
so that the MDPDE of $(\theta, \gamma)$ is to be obtained by minimizing 
$$
H_n(\theta, \gamma) = \iint V(y, x; \theta, \gamma)d\widehat{G}(x,y),
$$
with respect to the parameters. Further, the $\psi$-function for the MDPDEs 
as given by Equation (\ref{EQ:psi-2}) satisfies 
\begin{eqnarray}
	\psi_1(Y, X; \theta, \gamma)  = \frac{\partial V(Y, X; \theta, \gamma) }{\partial\theta}, 
	~~~~~~~~~~
	\psi_2(Y, X; \theta, \gamma) = \frac{\partial V(Y, X; \theta, \gamma) }{\partial\gamma}. 
\end{eqnarray}
Now, let us assume the following conditions:

\begin{itemize}
	\item[(D1)] The supports of the distributions $F_\theta$ and $F_{X,\gamma}$ for any value of $X$
	are independent of the parameters $\theta$ and $\gamma$ respectively.
	The true distribution $G(x,y)$ is also supported on the set $A = \{(x,y) : f_\theta(y|x)f_{X,\gamma}(x)>0 \}$,
	on which the true density $g$ is positive. 
	
	\item[(D2)] There exists an open subset $\omega$ of the parameter space 
	that contains the best  fitting parameter $(\theta_0, \gamma_0)$ and for all $(\theta, \gamma) \in \omega$ 
	and for almost all $(x,y) \in A$, the densities  $f_\theta$ and $f_{X,\gamma}$ are 
	thrice continuously differentiable with respect to $\theta$ and $\gamma$ respectively.
	
	\item[(D3)] The integrals $\iint f_\theta(y|x)^{1+\alpha}f_{X,\gamma}(x)^{1+\alpha}dxdy$ and
	$\iint f_\theta(y|x)^{\alpha}f_{X,\gamma}(x)^{\alpha}dG(x,y)$ can be differentiated three times 
	and the derivatives can be taken under the integral sign. 
	Further the $\psi$-function under consideration is finite.
	
	\item[(D4)] The matrix 
	\begin{eqnarray}
		\Lambda_G(\theta_0, \gamma_0) &=& \iint 
		\left.\frac{\partial}{\partial(\theta, \gamma)}\psi(y,x,\theta,\gamma)
		\right|_{(\theta, \gamma)=(\theta_0, \gamma_0)}dG(x,y) \nonumber\\
		&=& \iint \left.\frac{\partial^2}{\partial(\theta, \gamma)^2} 
		V(y,x,\theta,\gamma)\right|_{(\theta, \gamma)=(\theta_0, \gamma_0)}dG(x,y),\nonumber
	\end{eqnarray}
	exists finitely and is non-singular.
	
	\item[(D5)] For all  $(\theta, \gamma) \in \omega$, 
	each of the third derivatives of $V(y,x,\theta,\gamma)$ with respect to $(\theta, \gamma)$ 
	is bounded by a function of $(x,y)$, independent of $(\theta, \gamma)$, 
	that has finite expectation with respect to the true distribution $G$.
\end{itemize}

\begin{theorem}
	Under Assumptions (A1), (A2), (A5) and (D1)--(D5), 
	there exists a sequence of solutions  $\{(\hat{\theta}_n, \hat{\gamma}_n)\}$
	of the minimum density power divergence estimating equations (\ref{EQ:est-eqn-1}) and (\ref{EQ:est-eqn-2})
	with probability tending to one, that is consistent for the best fitting parameter $(\theta_0, \gamma_0)$.
	\\
	(Then, the asymptotic normality of this sequence $\{(\hat{\theta}_n, \hat{\gamma}_n)\}$ 
	follows from Theorem \ref{THM:AsymN_M-est} under the assumptions of that theorem.) 
	\label{THM:Asymp_MDPDE}
\end{theorem}
\noindent
\textbf{Proof: }
We follow a similar argument to that in the proof of Theorem 6.4.1(i) of Lehman (1983).
Consider the behavior of $H_n(\theta, \gamma)$, as a function of $(\theta, \gamma)$,
on a sphere $Q_a$ having center $(\theta_0, \gamma_0)$ and radius $a$. 
Then, to prove the existence part, it is enough to show that, for sufficiently small $a$, 
\begin{eqnarray}
	H_n(\theta, \gamma) > H_n(\theta_0, \gamma_0),
	\label{EQ:asymp_MDPDE_to_show}
\end{eqnarray} 
with probability tending to one, 
for any point  $(\theta, \gamma)$ on the surface of $Q_a$. 
Hence, for any $a>0$,  $H_n(\theta, \gamma)$ has a local minimum in the interior of $Q_a$ and
the estimating equations of the MDPDE have a solution  $\{(\hat{\theta}_n(a), \hat{\gamma}_n(a)\}$
within $Q_a$, with probability tending to one.

Now a Taylor series expansion of  $H_n(\theta, \gamma)$ around $(\theta_0, \gamma_0)$ yields
\begin{eqnarray}
	H_n(\theta_0, \gamma_0) - H_n(\theta, \gamma) 
	&=& - \sum_{i=1}^{q+r} (\zeta_i - \zeta_i^0)
	\left.\frac{\partial H_n(\theta,\gamma)}{\partial\zeta_i}\right|_{(\theta, \gamma)=(\theta_0, \gamma_0)}
	\nonumber\\
	&& - \frac{1}{2}\sum_{i, j =1}^{q+r} (\zeta_i - \zeta_i^0) (\zeta_j - \zeta_j^0)
	\left.\frac{\partial^2 H_n(\theta,\gamma)}{\partial\zeta_i\zeta_j}\right|_{(\theta, \gamma)=(\theta_0, \gamma_0)}
	\nonumber\\
	&& + \frac{1}{6}\sum_{i, j, k =1}^{q+r} (\zeta_i - \zeta_i^0) (\zeta_j - \zeta_j^0)(\zeta_k - \zeta_k^0)
	\left.\frac{\partial^3 H_n(\theta,\gamma)}{\partial\zeta_i\zeta_j\zeta_k}\right|_{(\theta,\gamma)=(\theta^*,\gamma^*)}
	\nonumber\\
	&=& S_1 + S_2+S_3, ~~~~ \mbox{(say),}
\end{eqnarray}
where $\zeta_i$ and $\zeta_i^0$ are the $i$-th component of the parameter vectors 
$(\theta, \gamma)$ and $(\theta_0, \gamma_0)$ respectively for all $i=1, \ldots, q+r$, 
and $(\theta^*,\gamma^*)$ lies in between $(\theta, \gamma)$ and $(\theta_0, \gamma_0)$ 
with respect to the Euclidean norm.
By a direct extension of the arguments presented in the proof of Theorem 3.1 of Basu et al.~(2006), 
we get, with probability tending to one, on $Q_a$,
\begin{eqnarray}
	|S_1| &<& (q+r) a^3, ~~\mbox{for all $a>0$}; ~~~~~~ \nonumber\\
	S_2 &<& -c a^2, ~~\mbox{for all $a<a_0$ with some $c, a_0 >0$}; ~~~~~~\nonumber\\
	\mbox{and}~~~|S_3| &<& b a^3, ~~\mbox{for all $a>0$ with some $b>0$},\nonumber
\end{eqnarray}
using the assumptions (D1)--(D5) and Lemma \ref{LEM:asym_lambda} whenever necessary.
Combining these, we get 
$$\max(S_1+S_2+S_3) < -ca^2 + (b+q+r)a^3,$$ 
which is less that zero whenever $a < \frac{c}{b+q+r}$ proving (\ref{EQ:asymp_MDPDE_to_show}) holds.

Finally, to show that one can choose a root of the estimating equations of MDPDEs independent of the 
radius $a$, consider the sequence of roots closest to the best fitting parameter $(\theta_0, \gamma_0)$,  
which exists by continuity of  $H_n(\theta, \gamma)$ as a function of $(\theta, \gamma)$. 
This sequence will also be consistent completing the proof of the theorem.
\hfill{$\square$}

\bigskip
Note that Assumptions (D1)--(D5) are easier to check compared to the
(stronger) Assumptions (A6)--(A9) and are the routine extensions of the 
corresponding assumptions [(A1)--(A5)] of Basu et al.~(2006).

\subsection{Other M-estimators with Different $\psi$-Functions}

Although our main focus in this paper is to study one particular M-estimator, 
namely the minimum density power divergence estimator (MDPDE), 
it opens the scope of many different M-estimators through the general 
results derived in Section \ref{SEC:SC_AN_M-est}. This general framework of parameter estimation
based on some suitable estimating equation is well studied in case of complete data
and several optimum robustness properties of these M-estimators has been proved 
for different classes of weight function; for example, see Huber (1981) and Hampel et al.~(1986).
In fact, there exists different class of $\psi$-function generating optimum solution in case of different problems.
For example, in case of estimating the location parameter in a symmetric distribution, 
the $\psi$ functions, that are  odd in the targeted parameter,
lead to such optimum M-estimation. 

However, as pointed out in Wang (1999), 
an optimum $\psi$-function for the complete data might not enjoy similar optimality for the 
censored data, even if there is no covariable presence.
The main reason is that the lifetime variables are not usually symmetric and neither belong 
to a location-scale family; rather it is usually asymmetric. 
The case of censored data with covariates, as considered here, is much more complicated and 
we can not directly pick a $\psi$-function from the theory of complete data.
Wang (1999) presented some example of $\psi$-functions in the context of censored data
with no covariates that can be extended in the present case with several covariables.
However, their usefulness and optimality both in terms of efficiency and robustness need to be verified
for the censored data cases with or without covariates. 
There need a lot of research in this area to suggest an optimum $\psi$-function 
under any suitable criteria of robustness or efficiency based on censored data.

However, we believe that the minimum density power divergence estimator
proposed here is quite sufficient for most practical situations
since it produces highly robust estimators with only a slight loss in efficiency compared to the 
maximum likelihood estimator (as described in Section \ref{SEC:Numerical}). 
Further, the estimating equation of MDPDEs can be solved by 
any simple numerical technique quite comfortably and has a simple interpretation in terms of the 
density power divergence.  Thus, although some future research work 
may provide suitable $\psi$-function satisfying some optimality criteria with complicated form or 
estimation procedure, the MDPDE will  still have its importance in many practical scenarios due to its simplicity.


\section{Global nature and the Semi-parametric Extensions}\label{SEC:SemiP_Ex}

Although the main focus of our paper is the MDPDE under fully parametric set-up, the M-estimator defined in Definition \ref{DEF:M-est}
and its asymptotic theory derived in the previous section is completely general in the sense that 
it can also be applied to any semi-parametric or even non-parametric set-ups. 
To see this, just note that the general M-estimator is defined in terms of a $\psi$ function 
that only need to satisfy Equation (\ref{EQ:psi-prop1}). Therefore, one can also consider the $\psi$ functions,
$\psi(y,x;\theta)$, involving no parametric assumptions on the distribution of $x$  (and hence independent of parameter $\gamma$)
and define the  M-estimator as before based on the corresponding estimating equation;
that estimator will also follow the general asymptotic theory developed in this paper.
Further, in this case, we might generalize our requirement (\ref{EQ:psi-prop1}) for such $\psi$ functions
by considering integral with respect to only the conditional distribution $G_{Y|X}$ of $Y$ given $X$ as follows 
(since $\phi$ doesn't include any distributional part of $X$): 
\begin{equation}
	\int  \psi(y, x; \theta) d{G_{Y|X}}(y) =0.
	\label{EQ:psi_prop2}
\end{equation}

In this general sense, the existing estimators of Zhou (2010) and Wang et al.~(2015)
become particular members of our class of general M-estimators with some specific choice of $\psi$ function 
without distributional assumptions on $X$. In particular the choice
\begin{eqnarray}
	\psi(y,x;\theta)= \psi_0\left(y-x^T\theta\right),
\end{eqnarray}
under the set-up considered in Section \ref{SEC:regression} (except the distributional assumption on $X$)
generates the estimator proposed in Zhou (2010).
Then the  asymptotic results  of Zhou et al (2010) directly follows from our general theory of Section \ref{SEC:M-estimators};
in particular, Theorem 3.3 of Zhaou (2010) follows from our Theorem 4.4.

Similarly, the proposal of Wang et al.~(2015) can also be though of as a special case of our general M-estimators
under the set-up of \ref{SEC:AFT} (except the distributional assumption on $X$) with the $\psi$ function
\begin{eqnarray}
	\psi(y,x;\theta) = \left[\begin{array}{c}
		\psi_0\left(\frac{\omega(x)(y-x^t\beta)}{\sigma}\right)x\omega(x)\\
		\chi\left(\frac{\omega(x)(y-x^t\beta)}{\sigma}\right)
	\end{array}\right],
\end{eqnarray}
where $\chi(s)=s\psi_0(s)-1$, $\omega(x)$ is some suitable weights and $\psi_0$ is some suitable function as given in Wang et al.~(2015).
Once again, all the asymptotic results of their paper follow from our general theory presented in Section \ref{SEC:M-estimators}.
For example, Theorem 3.2 and 3.3 of Wang et al.~(2015) follow from our Theorem 4.2 and 4.4. respectively
under the above mentioned set-up. 
A numerical comparison of our MDPDE with the estimator of Wang et al.~(2015) has been provided later in Section \ref{SEC:med_App}
through an interesting real data example.

However, the proposed MDPDE, a special M-estimator with the $\psi$ function given by (\ref{EQ:psi-2}),
involve the assumed density of the covariates $X$. So it cannot be applied directly to the 
semi-parametric settings where no distributional assumption has been made.
But, we can easily extend our definition of MDPDE for the semi-parametric cases by 
considering the density power divergence between the conditional densities $f_\theta(Y|X)$ of $Y$ given $X$ 
instead of considering the joint density of $Y$ and $X$. The $\psi$ function corresponding 
to this extended MDPDE under semi-parametric set-up can be seen to have the form
\begin{eqnarray}
	&=& \widetilde{\zeta}_\theta(X)  
	-  u_\theta(Y, X) f_\theta(Y|X)^{\alpha},  \label{EQ:psi-3}
\end{eqnarray}
where $\widetilde{\zeta}_\theta(x) = \int u_\theta(y, x) f_\theta(y|x)^{1+\alpha}dy$. 
Clearly, this $\psi$ function, $\psi_{\alpha}(y, x; \theta)$, corresponding to the extended MDPDE satisfies 
the stronger condition (\ref{EQ:psi_prop2}) and hence also satisfies (\ref{EQ:psi-prop1}). 
Thus, all the properties derived in Section \ref{SEC:M-estimators}
continue to hold under suitable modification for the semi-parametric set-up.
With this modification, the proposed MDPDE can now be applied to any semi-parametric set-up including the linear regression set-up
of Zhou et al.~(2010), as considered in Section \ref{SEC:regression} with fully parametric assumptions.

\section{Robustness: Influence Function Analysis}
\label{SEC:IF}

The influence function (Hampel et al., 1986) of an estimator is a popular tool to measure its classical robustness properties.
It measures the stability of the estimator under infinitesimal contamination yielding a first order approximation 
of the bias due to that small contamination in data. More precisely,  if $T_\psi(G)=(T_\psi^\theta(G), T_\psi^\gamma(G))$
denotes the statistical functional for the M-estimator corresponding to $\psi$ (which satisfies Equation (\ref{EQ:psi-prop1})),
then the influence function of this estimator is defined as
$$
IF((y_0,x_0);T_\psi, G) = \left.\frac{\partial}{\partial\epsilon}T_\psi(G_\epsilon)\right|_{\epsilon=0} 
=\lim\limits_{\epsilon\downarrow 0}\frac{T_\psi(G_\epsilon) - T_\psi(G)}{\epsilon},
$$
where $G_\epsilon=(1-\epsilon)G + \epsilon\wedge_{(x_0,y_0)}$ is the contaminated distribution
with $\epsilon$ being the contamination proportion and $\wedge_{(x_0,y_0)}$ being the degenerate distribution 
at the contamination point $(x_0,y_0)$. If the influence function is bounded in the contamination points $(x_0,y_0)$,
the bias under infinitesimal contamination cannot  become arbitrarily large even when the contamination is very far from the data center;
hence the estimator will be robust with respect to the  data contamination.

A straightforward albeit lengthy differentiation of the estimation equation 
(Equation (\ref{EQ:psi-prop1}) with $G$ replaced by $G_\epsilon$ and $(\theta,\gamma)$ replaced by $T_\psi(G_\epsilon)$) 
yields the form of the influence function of our general M-estimators, which is presented in the following theorem.

\begin{theorem}
	\label{THM:IF_Mest}
	Under the above mentioned set-up,
	\begin{eqnarray}
		IF((y_0,x_0);T_\psi, G) = \Lambda_G(T_\psi^\theta(G), T_\psi^\gamma(G))^{-1}\psi(y_0,x_0;T_\psi^\theta(G), T_\psi^\gamma(G)).
		\label{EQ:IF_Mest}
	\end{eqnarray}
\end{theorem}

Clearly, whenever we choose the $\psi$ function to be bounded with respect to $y$ and $x$,
the influence function of the corresponding M-estimator will be bounded in both $y_0$ and $x_0$;
hence the estimator will be robust with respect to both the outlier $y_0$ in response variables as well as 
the leverage point $x_0$ in the explanatory variables.
However, if we have $\psi$ function bounded only in $y$ and not in $x$ 
(like the $\psi$ functions of the classical M-estimators under normal linear regression without censoring)
the resulting estimator will be robust only with respect to outliers in response 
but may not be robust with respect to leverage points.

In particular, the above theorem also provides the influence function of the proposed MDPDE 
under fully parametric models by just using the $\psi$ function given in (\ref{EQ:psi-2}). Note that, for most common parametric models,
this particular $\psi$ function is bounded in both $y$ and $x$ whenever $\alpha>0$ implying the robust nature of 
the MDPDE with $\alpha>0$ for both the outliers in responses and covariates.
However, at $\alpha=0$, the $\psi$ function of the corresponding MDPDE (which is the same as the MLE) 
is proportional to the score functions which are generally unbounded for most parametric models 
and prove their non-robust nature.

However, the MDPDE under the semi-parametric extension has a different  $\psi$ function, given in (\ref{EQ:psi-3}), 
which is not bounded in $x$ for all $\alpha\geq0$ under common parametric families;
but it is generally bounded in $y$ for $\alpha>0$. Hence the semi-parametric MDPDE with $\alpha>0$ are 
only robust with respect to outliers in the response but not robust under leverage points.

%
%

\section{Numerical Illustrations}
\label{SEC:Numerical}

\subsection{Simulation Study}

Consider the exponential regression model with randomly censored data 
and normal covariables as discussed in the previous subsection. For simulation exercise,
we consider only one covariable so that $X$  is a univariate normal random variable with 
mean $\gamma$ (scalar) and variance $1$. Then a covariate sample of size $n$ is generated from 
$N(\gamma,1)$ distribution and given the value $x$ of the covariate we simulate the (lifetime) response variable 
from an exponential distribution with mean $\theta x$ under a random censoring scheme; 
the true values of the parameters are taken to be $\theta=1$ and $\gamma=5$.  
Here we consider the simple exponential censoring distribution, 
but the censoring rate is determined to keep the expected proportion of censoring at 10 or $20\%$ 
under the true distribution.
Under the exponential censoring distribution with mean $\tau$, i.e., $C \sim Exp(\tau)$,
the expected proportion of censoring under the true distribution $Exp(\theta x)$ can be seen to be
$$
P(Y>C) = \frac{\theta x}{\tau + \theta x}.
$$
So to make this proportion equal to $10\%$ or $20\%$, 
we need to take $ \tau = 9\theta x$ and $\tau = 4\theta x$ respectively
(with $\theta=1$ for our simulation study).

Then we compute the MDPDE of $(\theta, \gamma)$ numerically and repeat the process  
$1000$ times to obtain the empirical estimates for the total absolute bias 
(sum of the absolute biases of $\theta$ and $\gamma$) and the total MSE
(sum of the MSEs of $\theta$ and $\gamma$) of the MDPDE with respect to the target value $(1, 5)$.
One can measure the performance of the proposed MDPDE with respect to several choice of tuning parameter $\alpha$ 
and with the maximum likelihood estimator (MLE) at $\alpha=0$ by comparing these empirical bias and MSEs.

At first we consider only the pure sample without any contamination and compare the
efficiencies of the MDPDEs for different $\alpha$ with the MLE (at $\alpha=0$). 
The empirical estimates of efficiency are computed from the total MSEs and are reported in Table \ref{TAB:MDPDE_noCont}
along with the total absolute bias for different $\alpha$ and different censoring proportions.
It is clear from the table that the efficiency of the MDPDE decreases as $\alpha$ increases
but the loss in efficiency is not so significant at smaller positive values of $\alpha$.
Further, for any fixed $\alpha$ both the total absolute bias and MSE increase 
as the censoring proposing increases.

\begin{table}[h]
	\caption{Empirical Summary measures for the MDPDEs  under no contamination}
	\resizebox{\textwidth}{!}{
		\begin{tabular}{l|l|rrrrrrr}\hline
			& & \multicolumn{7}{c}{$\alpha$}\\
			&	Cens. Prop.	&	0.00	&	0.01	&	0.10	&	0.30	&	0.50	&	0.70	&	1.00	\\\hline
			Total Abs.	&	10\%	&	0.3848 & 0.3695 & 0.4067 & 0.4659 & 0.5145 & 0.5474 & 0.5801 	\\
			Bias			&	20\%	&	0.4260	&	0.4192	&	0.4616	&	0.5472	&	0.6122	&	0.6518	&	0.6890	\\
			\hline
			Total			&	10\%	&	0.1363 & 0.1490 & 0.1577 & 0.1873 & 0.2198 & 0.2458 & 0.2773	\\
			MSE			   &	20\%	&	0.1860	&	0.1954	&	0.2078	&	0.2656	&	0.3097	&	0.3424	&	0.3782	\\
			\hline
			Relative	  &		10\%	&	100\%	&	91\% 	& 	86\% 	& 	73\% 	& 	62\% 	& 	55\% 	& 	49\% \\
			Efficiency	  &		20\%	&	100\%	&	95\%	&	90\%	&	70\%	&	60\%	&	54\%	&	49\%	\\
			\hline
		\end{tabular}}
		\label{TAB:MDPDE_noCont}
	\end{table}

	Next, to examine the robustness of the proposed MDPDEs over the MLE, we repeat the
	above simulation study but with 5, 10, 15 or $20\%$ contamination in the response variable and covariates.
	For contamination in response variable, we generate them from an $Exp(5x)$ distribution ($\theta=5$)
	under the same censoring scheme as before;
	for contamination in the covariates, we simulate observations from another normal distribution 
	with mean $10$ and variance $1$. The empirical bias and MSE of the estimators are reported in
	Tables \ref{TAB:MDPDE_Bias_Cont} and \ref{TAB:MDPDE_MSE_Cont} respectively.
	Clearly, note that the total absolute bias as well as the total  MSE increases for any fixed $\alpha$
	as the contamination proportion increases. However, these changes are rather drastic 
	at smaller values of $\alpha$ and stabilize as $\alpha$ increases. 
	In other words, the MDPDE with larger $\alpha \geq 0.3$ can successfully ignore the outliers 
	to generate robust inference.

	\begin{table}[h]
		\caption{Empirical total absolute bias of the MDPDEs  for different contamination proportions}
		\resizebox{\textwidth}{!}{
			\begin{tabular}{l|l|rrrrrrr}\hline
				& & \multicolumn{7}{c}{$\alpha$}\\
				Cens. Prop.	&	Cont. Prop	&	0.00	&	0.01	&	0.10	&	0.30	&	0.50	&	0.70	&	1.00	\\\hline
				10\%	&	5\%	&	0.513	&	0.506	&	0.392	&	0.270	&	0.220	&	0.199	&	0.186	\\
				&	10\%	&	0.741	&	0.667	&	0.437	&	0.258	&	0.223	&	0.216	&	0.217	\\
				&	15\%	&	1.093 & 1.011 & 0.794 & 0.554 & 0.518 & 0.534 & 0.578	\\
				&	20\%	&	1.594	&	1.476	&	1.169	&	1.009	&	0.878	&	0.816	&	0.781	\\\hline
				20\%	&	5\%	&	1.006	&	0.953	&	0.747	&	0.651	&	0.723	&	0.786	&	0.856	\\
				&	10\%	&	0.759	&	0.706	&	0.522	&	0.412	&	0.393	&	0.405	&	0.415	\\
				&	15\%	&	0.865	&	0.790	&	0.719	&	0.559	&	0.481	&	0.449	&	0.436	\\
				&	20\%	&	1.090	&	1.000	&	1.102	&	1.166	&	1.092	&	1.035	&	1.003	\\
				\hline
			\end{tabular}}
			\label{TAB:MDPDE_Bias_Cont}
		\end{table}
		
		\begin{table}[h]
			\caption{Empirical total MSE of the MDPDEs  for different contamination proportions}
			\resizebox{\textwidth}{!}{
				\begin{tabular}{l|l|rrrrrrr}\hline
					& & \multicolumn{7}{c}{$\alpha$}\\
					Cens. Prop.	&	Cont. Prop	&	0.00	&	0.01	&	0.10	&	0.30	&	0.50	&	0.70	&	1.00	\\\hline
					10\%	&	5\%	&	0.265	&	0.257	&	0.139	&	0.081	&	0.070	&	0.071	&	0.083	\\
					&	10\%	&	0.453	&	0.425	&	0.220	&	0.101	&	0.079	&	0.077	&	0.081	\\
					&	15\%	&	0.908	& 0.848 & 0.497 & 0.259 & 0.229 & 0.244 & 0.298	\\
					&	20\%	&	2.282	&	2.167	&	1.585	&	0.913	&	0.714	&	0.629	&	0.602	\\\hline
					20\%	&	5\%	&	0.794	&	0.743	&	0.425	&	0.305	&	0.353	&	0.410	&	0.482	\\
					&	10\%	&	0.441	&	0.415	&	0.248	&	0.166	&	0.144	&	0.143	&	0.147	\\
					&	15\%	&	0.689	&	0.680	&	0.499	&	0.309	&	0.249	&	0.235	&	0.242	\\
					&	20\%	&	1.208	&	1.171	&	1.153	&	1.008	&	0.898	&	0.838	&	0.833	\\
					\hline
				\end{tabular}}
				\label{TAB:MDPDE_MSE_Cont}
			\end{table}

			\subsection{Real Data Application : Heart Transplant data}
			\label{SEC:med_App}
			
			We will now apply our proposed MDPDEs to an interesting real data example with the semi-parametric model assumptions,
			which will illustrate the performance of the proposed semi-parametric  extension described in Section \ref{SEC:SemiP_Ex}
			along with its applicability in real life scenarios.
			The data set considered is from the popular Stanford heart transplant program described in details in Clark et al.~(1971) 
			and contains the following survival information of 158 patients (Crowley and Hu, 1977; Escobar and Meeker Jr, 1992):
			ID number of patients (``ID"), survival or censoring time (``TIME"), censoring status (dead or alive), 
			patient's age at first transplant in years (``AGE") and the T5 mismatch score (``T5-MS").
			The dataset, available from the `survival' library of R, 
			was analyzed statistically by many authors including Brown et al.~(1973), Turnbull et al.~(1974), Mantel and Byar (1974), 
			and Miller and Halpern (1982).  Recently, it has also been used to illustrate the performances of robust estimates under 
			semi-parametric AFT models by Salibian-Barrera and Yohai (2008), Locatelli et al.~(2011) and Wang et al.~(2015).
			The latest robust estimator of Wang et al.~(2015), namely, the KMW-GM estimator, has been seen to work best for this dataset 
			while using the model 
			$$
			\log({\rm TIME}) = \beta_0 + \beta_1 ({\rm AGE}) + \beta_2 (\mbox{T5-MS}) + \sigma\epsilon. 
			$$
			Here, we will apply our proposed MDPDE with different tuning parameters $\alpha$ with the same parametric model as above 
			and the assumption that $\epsilon\sim N(0,1)$ and illustrate the superior performance of our proposal 
			over the KMW-GM estimator of Wang et al.~(2015).

			As noted in Wang et al.~(2015), there are three potential outliers in the dataset corresponding to the ID 2, 16 and 21, 
			where the patients have unexpectedly shorter survival times. 
			This finding is also consistent with the results from previous analyses of the dataset
			and so we also treat these three data points as outliers and compute our MDPDEs twice; 
			once with the full data set and once after removing these outliers. 
			However, since the estimates of the parameters $(\beta_0,~\beta_1,~\beta_2,~\sigma)$ differ only slightly 
			in the two cases with and without outliers, we will report teh relative variation in the estimates in order 
			to check the extent of their robustness. Following Wang et al.~(2015), we define the relative variation as
			$$
			{\rm Relative Variation} = \frac{|\hat{\theta}_{\rm full}-\hat{\theta}_{\rm cleaned}|}{|\hat{\theta}_{\rm full}|},
			$$
			where $\hat{\theta}_{\rm full}$ is the estimated parameter value based on the full data set and
			$\hat{\theta}_{\rm cleaned}$ is the parameter estimate based on the cleaned data after removing the three outliers.
			The relative variations obtained for each of the parameters are reported in Table \ref{TAB:RA_data_example}
			for our proposed MDPDEs with different tuning parameters and $\alpha$ and also for the KMW-GM estimator of Wang et al.~(2015).   
			It can be seen clearly from the table that the MDPDEs of most of the parameters are much more stable and have less relative variation
			compared to the KMW-GM estimator of Wang et al.~(2015) for $\alpha\geq 0.4$; this clearly shows the greater robustness of our proposal
			compared to the existing robust method. Further, note that the relative variation is quite high at $\alpha=0$ which is the 
			non-robust maximum likelihood estimator. As $\alpha$ increases the relative variations of all the parameters decrease
			significantly which again shows the significant gain in robustness of our proposal with increasing $\alpha$.

			\begin{table}[h]
				\caption{The relative variation of the MDPDEs at $\alpha$ and the KMW-GM estimates of Wang et al.~(2015) 
					with and without outliers for the Heart Transplant Data}
				\centering 
				\begin{tabular}{l|rrrr}\hline
					$\alpha$	&	$\beta_0$	&	$\beta_1$	&	$\beta_2$	&	$\sigma$	\\\hline
					0	&	0.0191	&	0.8198	&	0.3622	&	0.0886	\\
					0.05	&	0.0200	&	0.3450	&	0.3189	&	0.0786	\\
					0.1	&	0.0202	&	0.1477	&	0.2727	&	0.0707	\\
					0.2	&	0.0191	&	0.0041	&	0.1825	&	0.0570	\\
					0.3	&	0.0172	&	0.0347	&	0.1079	&	0.0463	\\
					0.4	&	0.0153	&	0.0445	&	0.0558	&	0.0384	\\
					0.5	&	0.0139	&	0.0452	&	0.0236	&	0.0328	\\
					0.7	&	0.0118	&	0.0409	&	0.0039	&	0.0259	\\
					0.9	&	0.0108	&	0.0369	&	0.0088	&	0.0225	\\
					1	&	0.0106	&	0.0354	&	0.0077	&	0.0215	\\\hline
					KMW-GM	&	0.0153	&	0.0034	&	0.0807	&	0.0573	\\
					\hline
				\end{tabular}
				\label{TAB:RA_data_example}
			\end{table}

			\section{On the Choice of Tunning Parameter $\alpha$ in MDPDE}\label{SEC:selection_alpha}
			
			A crucial issue for applying the proposed MDPDE in any real-life problem
			is the choice of tuning parameter $\alpha$. As we have seen that the
			robustness and efficiency of the MDPDEs depend crucially on the tuning parameter $\alpha$, 
			it needs to be chosen carefully in practice where we have no idea regarding the contamination 
			and censoring proportions. The simulation study presented in Section \ref{SEC:Numerical} gives some 
			indication in this direction. We have seen that the MDPDEs with larger $\alpha \geq 0.3$ 
			are robust enough to successfully address the problem of outliers; the robustness increases as $\alpha$ increases.
			On the other hand, the efficiency of the MDPDEs under pure data is seen to decreases as $\alpha$ increases,
			but there is no significant loss in efficiency at smaller positive values of $\alpha$ near 0.3.
			So, we recommend to use a value of the tuning parameter $\alpha$ near 0.3 
			to get a fair compromise between efficiency and robustness whenever the amount of contamination is not known in practice.
			This is in-line with the  empirical suggestions given by Basu et al.~(2006) in the context of MDPDE
			based on censored data with no covariables. 
			However, these empirical suggestions need further justification based on more elaborative simulation 
			and theoretical aspects. In case of complete data, some such justifications of the 
			data driven choice of $\alpha$ is given by Hong and Kim (2001) and Warwick and Jones (2005). 
			Their work might have been generalized to the case of censored data, although it is not very easy, 
			in order to solve this issue of selecting $\alpha$. We hope to pursue this in our future research.

			\section{Conclusion}\label{SEC:Conclusion}
			
			The present paper proposes the minimum density power divergence estimator 
			under the parametric set-up for censored data with covariables to generate 
			highly efficient and robust inference. The applicability of the proposed technique 
			is illustrated through appropriate theoretical results and simulation exercise in 
			the context of censored regression with  stochastic covariates. 
			Further, the paper provide the asymptotic theory for a general class of estimators
			based on the estimating equation which opens the scope of studying many such estimators
			in the context of censored data in presence of some stochastic covariates.

			\bigskip\noindent
			\textbf{Acknowledgments:} 
			The authors gratefully acknowledge the comments of two anonymous referees which led to an improved version of the manuscript.



\begin{thebibliography}{}
				
				\bibitem{}
				Basu, S., Basu, A., and Jones, M. C. (2006). Robust and efficient parametric estimation for censored survival data. 
				{\em Annals of the Institute of Statistical Mathematics}, {\bf 58(2)}, 341--355.
				
				\bibitem{}
				Basu, A., Harris, I. R., Hjort, N. L., and Jones, M. C. (1998). 
				Robust and efficient estimation by minimising a density power divergence. 
				{\em Biometrika}, {\bf 85(3)}, 549--559.
				
				\bibitem{}
				Bednarski, T. (1993)
				Robust estimation in the Cox regression model.
				{\em Scand. J. Statist.}, {\bf 20}, 213--225.
				
				\bibitem{}
				Bednarski, T. and Borowicz, F. (2006). 
				coxrobust: Robust Estimation in Cox Model. {\em R package version 1.0}.
				
				\bibitem{}
				Begun, J. M., Hall, W. J., Huang, W. M., Wellner, J. A. (1983).
				Information and Asymptotic Efficiency in Parametric-Nonparametric Models.
				{\em Annals of Statistics}, {\bf 11}, 432--452.
				
				
				\bibitem{}
				Brown, B. W., Jr., Hollander, M., and Korwar, R. M. (1973). 
				Nonparametric Test of independence for censored data with application to Heart Transplant Studies. 
				{\em Florida State University Conference on Reliability and Biometry}. 
				
				
				
				\bibitem{}
				Buckley, J., and James, I., (1979). 
				Linear regression with censored data. 
				{\em Biometrika}, {\bf 66}, 429--436.
				
				\bibitem{}
				Cai, Z. (1998). Asymptotic properties of Kaplan-Meier estimator for censored dependent data. 
				{\it Statistics and probability letters}, {\bf 37(4)}, 381--389.
				
				\bibitem{}
				Campbell, G., and F\"{o}ldes, A. (1982). 
				Large sample properties of nonparametric bivariate estimators with censored data. 
				{\it Nonparametric statistical inference}, {\bf 1}, 103-121.
				
				
				\bibitem{}
				Chen, Y. Y., Hollander, M., and Langberg, N. A. (1982). 
				Small-sample results for the Kaplan-Meier estimator. 
				{\em Journal of the American Statistical Association}, {\bf 77}, 141-144.
				
				
				\bibitem{}
				Clark, D. A., Stinson, E. B., Griepp, R. B., Schroeder, J. S., Shumway, N. E.,  and Harrison, D. C. (1971). 
				Cardiac Transplantation in Man. VI. Prognosis of Patients Selected  for Cardiac Transplantation.
				{\em Annals of Internal Medicine}, {\bf 75}, 15--21. 
				
				
				
				\bibitem{CO:84}
				Collett, D.  (2003).
				{\em Modelling Survival Data in Medical Research}.
				Chapman Hall, London, U.K.
				
				\bibitem{}
				Cox, D.R. (1972). Regression models and life tables (with discussion). 
				{\em Journal of Royal Statistical Society, Series B}. {\bf 34}, 187--220.
				
				
				
				
				\bibitem{CO:84}
				Cox, D. R., and Oakes, D.  (1984).
				{\em Analysis of Survival Data}.
				Chapman Hall, London, U.K.
				
				
				\bibitem{CO:91}
				Crowder, M. J., Kimber, A. C., Smith, R. L., and Sweeting, T. J. (1991).
				{\em Statistical Analysis of Reliability Data}.
				Chapman Hall, London, U.K.
				
				\bibitem{}
				Crowley, J. and Hu, M. (1977). 
				Covariance analysis of heart transplant survival data. 
				{\em Journal of the American Statistical Association}, {\bf 72}, 27--36. 
				
				
				\bibitem{}
				Dabrowska, D. M. (1988). Kaplan-Meier estimate on the plane. 
				{\em The Annals of Statistics}, {\bf 16(4)}, 1475--1489.
				
				\bibitem{}
				Escobar, L. A. and Meeker Jr, W. Q. (1992). 
				Assessing influence in regression analysis with censored data. 
				{\em Biometrics}, {\bf 48}, 507--528.
				
				\bibitem{}
				Farcomeni, A. and Viviani, S. (2011)
				Robust estimation for the Cox regression model based on trimming.
				{\em Biometrical Journal}, {\bf 53(6)}, 956--973.
				
				
				\bibitem{}
				Ghosh, A., and Basu, A. (2013). 
				Robust estimation for independent non-homogeneous observations using density power divergence 
				with applications to linear regression. {\em Electronic Journal of statistics}, {\bf 7}, 2420-2456.
				
				\bibitem{}
				Ghosh, A., and Basu, A. (2014). 
				Robust Estimation in Generalized Linear Models : The Density Power Divergence Approach. 
				{\em Test}, doi:10.1007/s11749-015-0445-3.
				
				
				\bibitem{Hampel/etc:1986}
				Hampel, F.~R., E.~Ronchetti, P.~J. Rousseeuw, and W.~Stahel (1986).
				{\em Robust Statistics: The Approach Based on Influence Functions}.
				New York, USA: John Wiley \& Sons.
				
				
				\bibitem{hk01}
				Hong, C. and Kim, Y. (2001), 
				Automatic selection of the tuning parameter in the minimum density power divergence estimation. 
				{\it Journal of the Korean Statistical Society}, {\bf 30}, 453--465.
				
				
				
				\bibitem{hlm:08}
				Hosmer, D. W., Lemeshow, S. and May, S.  (2008).
				{\em Applied Survival Analysis: Regression Modeling of Time-to-Event Data}.
				John Wiley \& Sons.
				
				
				\bibitem{Huber:1981}
				Huber, P.~J. (1981).
				{\em Robust Statistics}.
				John Wiley \& Sons.
				
				
				\bibitem{}
				Kaplan, E.~L., and Meier, P.~(1958). 
				Nonparametric estimation from incomplete observations. 
				{\it Journal of the American statistical association}, {\bf 53 (282)}, 457--481.
				
				
				
				\bibitem{}
				Kim, M., and Lee, S. (2008). Estimation of a tail index based on minimum density power divergence. 
				{\em Journal of Multivariate Analysis}, {\bf 99(10)}, 2453--2471.
				
				\bibitem{}
				Klein, J.P. and Moeschberger, M.L. (2003). 
				{\em Survival Analysis Techniques for Censored and Truncated Data, Second Edition.} 
				Springer-Verlag, New York.
				
				\bibitem{}
				Kosorok, M.R., Lee, B.L. and Fine, J.P. (2004). Robust inference for univariate proportional
				hazards frailty regression models. {\em Annals of Statistics}, {\bf 32}, 1448--1491.
				
				
				\bibitem{}
				Lawless, J.F. (2003). 
				{\em Statistical Models and Methods for Lifetime Data, Second Edition.} 
				John Wiley \& Sons, Inc. New York.
				
				
				\bibitem{}
				Lee, S., and Song, J. (2009). Minimum density power divergence estimator for GARCH models. 
				{\em Test}, {\bf 18(2)}, 316--341.
				
				
				\bibitem{}
				Lee, S., and Song, J. (2013). Minimum density power divergence estimator for diffusion processes. 
				{\em Annals of the Institute of Statistical Mathematics}, {\bf 65(2)}, 213-236.
				
				
				\bibitem{Lehmann:1983}
				Lehmann, E.~L. (1983).
				{\em Theory of Point Estimation}.
				John Wiley \& Sons.
				
				\bibitem{}
				Lo, S. H., Mack, Y. P., and Wang, J. L. (1989). 
				Density and hazard rate estimation for censored data via strong representation of the Kaplan-Meier estimator. 
				{\em Probability theory and related fields}, {\bf 80(3)}, 461-473.
				
				
				\bibitem{}
				Locatelli, I., Marazzi, A., Yohai, V. J. (2011). 
				Robust accelerated failure time regression.
				{\em Computational Statistics and Data Analysis}. {\bf 55}, 874--887.
				
				\bibitem{}
				Mantel, N. and Byar, D. P. (1974). 
				Evaluation of Response-Time data involving transient states: An illustration using Heart-Transplant data. 
				{\em Journal of the American Statistical Association}, {\bf 69}, 81--86. 
				
				
				
				\bibitem{}
				Miller, R. and Halpern, J. (1982). 
				Regression with censored data. 
				{\em Biometrika}, {\bf 69}, 521--531.
				
				\bibitem{}
				Peterson Jr, A.~V.~(1977). Expressing the Kaplan-Meier estimator as a function of empirical subsurvival functions. 
				{\it Journal of the American Statistical Association}, {\bf 72 (360a)}, 854--858.
				
				
				\bibitem{}
				Ritov, Y. (1986). Estimation in a Linear Regression Model with Censored Data. 
				{\em The Annals of Statistics}, {\bf 18(1)}, 303--328.
				
				\bibitem{}
				Robins, J. M., and Rotnitzky, A. (1992). 
				Recovery of information and adjustment for dependent censoring using surrogate markers. 
				In {\em AIDS Epidemiology}, 297--331. Birkh\''{a}user Boston.
				
				
				\bibitem{}
				Satten, G. A., and Datta, S. (2001). 
				The Kaplan–Meier estimator as an inverse-probability-of-censoring weighted average. 
				{\em The American Statistician}, {\bf 55(3)}, 207--210.
				
				
				\bibitem{}
				Salibian-Barrera, M., and Yohai, V. J. (2008). 
				High breakdown point robust regression with censored data. 
				{\em The Annals of Statistics}, {\bf 36(1)}, 118--146.
				
				
				\bibitem{Serfling:1980}
				Serfling, R.~J. (1980).
				{\em Approximation Theorems of Mathematical Statistics}.
				New York, USA: John Wiley \& Sons.
				
				
				\bibitem{}
				Stute, W. (1993). Consistent estimation under random censorship when covariables are present. 
				{\em Journal of Multivariate Analysis}, {\bf 45(1)}, 89--103.
				
				\bibitem{}
				Stute, W. (1995). The central limit theorem under random censorship. 
				{\em The Annals of Statistics}, {\bf 23(2)}, 422--439.
				
				\bibitem{}
				Stute, W. (1996). Distributional convergence under random censorship when covariables are present. 
				{\em Scandinavian Journal of Statistics}, {\bf 23(4)}, 461--471.
				
				\bibitem{}
				Stute, W., and Wang, J. L. (1993). The strong law under random censorship. 
				{\em The Annals of Statistics}, {\bf 21(3)}, 1591--1607.
				
				
				\bibitem{}
				Turnbull, B. W., Brown, B. W., Jr.,  and Hu, M. (1974). 
				Survivorship analysis of Heart Transplant data. 
				{\em Journal of the American Statistical Association}, {\bf 69}, 74--80. 
				
				
				\bibitem{}
				Tsai, W. Y., Jewell, N. P., and Wang, M. C. (1987). 
				A note on the product-limit estimator under right censoring and left truncation. 
				{\em Biometrika}, {\bf 74(4)}, 883--886.
				
				
				\bibitem{}
				Van der Laan, M. J., and Robins, J. M. (2003). 
				{\em Unified methods for censored longitudinal data and causality}. Springer.
				
				\bibitem{}
				Wang, J. L. (1999). Asymptotic Properties of M-Estimators Based on Estimating Equations and Censored Data. 
				{\em Scandinavian journal of statistics}, {\bf 26(2)}, 297--318.
				
				\bibitem{}
				Wang, M. C., Jewell, N. P., and Tsai, W. Y. (1986). 
				Asymptotic properties of the product limit estimate under random truncation. 
				{\em The Annals of Statistics}, {\bf 14(4)}, 1597--1605.
				
				
				\bibitem{}
				Warwick, J., and Jones, M. C. (2005). 
				Choosing a robustness tuning parameter. 
				{\em Journal of Statistical Computation and Simulation}, {\bf 75(7)}, 581--588.
				
				
				\bibitem{}
				Zhou, M. (1991). 
				Some properties of the Kaplan-Meier estimator for independent nonidentically distributed random variables. 
				{\em The Annals of Statistics}, {\bf 19(4)}, 2266--2274.
				
				
				\bibitem{}
				Zhou, M. (1992). M-estimation in censored linear models. Biometrika, 79(4), 837-841.
				
				
			\end{thebibliography}
\end{document}